# LOCAL WHITTLE ESTIMATION IN NONSTATIONARY AND UNIT ROOT CASES

By Peter C. B. Phillips[1] and Katsumi Shimotsu[2]

*Yale University, University of Auckland and University of York,*
*and Queen's University*

Asymptotic properties of the local Whittle estimator in the nonstationary case ($d > \frac{1}{2}$) are explored. For $\frac{1}{2} < d \leq 1$, the estimator is shown to be consistent, and its limit distribution and the rate of convergence depend on the value of $d$. For $d = 1$, the limit distribution is mixed normal. For $d > 1$ and when the process has a polynomial trend of order $\alpha > \frac{1}{2}$, the estimator is shown to be inconsistent and to converge in probability to unity.

**1. Introduction.** Semiparametric estimation of the memory parameter ($d$) in fractionally integrated $[I(d)]$ time series has attracted much recent study and is attractive in empirical applications because of its general treatment of the short memory component. Two commonly used semiparametric estimators are log periodogram (LP) regression and local Whittle estimation. LP regression is popular mainly because of the simplicity of its construction as a linear regression estimator. Local Whittle estimation involves numerical methods but is more efficient than LP regression. The local Whittle estimator was proposed by Künsch (1987) and Robinson (1995) showed its consistency and asymptotic normality for $d \in (-\frac{1}{2}, \frac{1}{2})$. Velasco (1999) extended Robinson's results to show that the estimator is consistent for $d \in (-\frac{1}{2}, 1)$ and asymptotically normally distributed for $d \in (-\frac{1}{2}, \frac{3}{4})$.

The present paper studies the asymptotic properties of the local Whittle estimator in the nonstationary case for $d > \frac{1}{2}$, including the unit root case and the case where the process has a polynomial time trend. These cases are of high importance in empirical work especially with economic time series,

Received July 2000; revised June 2003.

[1] Supported by NSF Grants SBR-97-30295 and SES 00-92509.

[2] Supported by a Cowles Foundation Cowles prize and by a Sloan Foundation fellowship.

which commonly exhibit nonstationary behavior and show some evidence of deterministic trends as well as long range dependence. The asymptotic properties of the local Whittle estimator in the nonstationary case over the region $d \in (\frac{1}{2}, 1)$ were explored in Velasco (1999). Velasco also showed that, upon adequate tapering of the observations, the region of consistent estimation of $d$ may be extended but with corresponding increases in the variance of the limit distribution. For the region $d \geq 1$, there is presently no theory for the untapered Whittle estimator and, for the region $d \in (\frac{3}{4}, 1)$, no limit distribution theory. The unit root case is of particular interest because it stands as an important special case of an $I(d)$ process with $d = 1$ and it has played a central role in the study of nonstationary economic time series. It is also now known to be the borderline that separates cases of consistent and inconsistent estimation by LP regression [Kim and Phillips (1999)] and, as we shall show here, local Whittle estimation.

This paper demonstrates that the local Whittle estimator (i) is consistent for $d \in (\frac{1}{2}, 1]$, (ii) is asymptotically normally distributed for $d \in (\frac{1}{2}, \frac{3}{4})$, (iii) has a non-normal limit distribution for $d \in [\frac{3}{4}, 1)$, (iv) has a mixed normal limit distribution for $d = 1$, (v) converges to unity in probability for $d > 1$ and (vi) converges to unity in probability when the process has a polynomial time trend of order $\alpha > \frac{1}{2}$. The present paper, therefore, complements the earlier work of Robinson (1995) and Velasco (1999) and largely completes the study of the asymptotic properties of the local Whittle estimator for regions of $d$ that are empirically relevant in most applications. The paper also serves as a counterpart to Phillips (1999b) and Kim and Phillips (1999), which analyze the asymptotics of LP regression for $d \in (\frac{1}{2}, 2)$.

The approach in the present paper draws on an exact representation and approximation theory for the discrete Fourier transform (d.f.t.) of nonstationary fractionally integrated processes. The theory, developed by Phillips (1999a), employs a model for nonstationary fractionally integrated processes that is valid for all values of $d$ and provides a uniform apparatus for analyzing the asymptotic behavior of their d.f.t.'s.

The remainder of the paper is organized as follows. Section 2 introduces the model. Consistency of the local Whittle estimator for $d \in (\frac{1}{2}, 1]$ and its inconsistency for $d > 1$ are demonstrated in Section 3. Section 4 derives the limit distributions. Results for fractionally integrated processes with a polynomial time trend are given in Section 5. Section 6 reports some simulation results and gives an empirical application using economic data. Section 7 makes some brief remarks on the important practical issue of finding a good general purpose estimator of $d$ when nonstationarity in the data is suspected. Some technical results are collected in Appendix A. Proofs are given in Appendix B.

**2. Preliminaries.** We consider the fractional process $X_t$ generated by the model

$$(1) \qquad (1-L)^d(X_t - X_0) = u_t, \qquad t = 0, 1, 2, \ldots,$$

where $X_0$ is a random variable with a certain fixed distribution. Our interest is in the case where $X_t$ is nonstationary and $d > \frac{1}{2}$, so in (1) we work from a given initial date $t = 0$, set $u_t = 0$ for all $t \leq 0$ and assume that $u_t$, $t \geq 1$, is stationary with zero mean and spectral density $f_u(\lambda)$. Expanding the binomial in (1) gives the form

$$(2) \qquad \sum_{k=0}^{t} \frac{(-d)_k}{k!}(X_{t-k} - X_0) = u_t,$$

where

$$(d)_k = \frac{\Gamma(d+k)}{\Gamma(d)} = (d)(d+1)\cdots(d+k-1)$$

is Pochhammer's symbol for the forward factorial function and $\Gamma(\cdot)$ is the gamma function. When $d$ is a positive integer, the series in (2) terminates, giving the usual formulae for the model (1) in terms of the differences and higher order differences of $X_t$. An alternate form for $X_t$ is obtained by inversion of (1), giving a valid representation for all values of $d$,

$$(3) \qquad X_t = (1-L)^{-d}u_t + X_0 = \sum_{k=0}^{t-1} \frac{(d)_k}{k!} u_{t-k} + X_0.$$

Define the discrete Fourier transform and the periodogram of a time series $a_t$ evaluated at the fundamental frequencies as

$$(4) \qquad w_a(\lambda_s) = \frac{1}{\sqrt{2\pi n}} \sum_{t=1}^{n} a_t e^{it\lambda_s}, \qquad \lambda_s = \frac{2\pi s}{n}, s = 1, \ldots, n,$$

$$I_a(\lambda_s) = |w_a(\lambda_s)|^2.$$

The model (1) is not the only model of nonstationary fractional integration. Another model that is used in the literature forms a process $X_t$ with $d \in [\frac{1}{2}, \frac{3}{2})$ from the partial sum of a stationary long-range dependent process, as in

$$(5) \qquad X_t = \sum_{k=1}^{t} U_k + X_0, \qquad d \in [\tfrac{1}{2}, \tfrac{3}{2}),$$

where $U_t$ has spectral density $f(\lambda) \sim G_0 \lambda^{-2(d-1)}$ as $\lambda \to 0$. Model (5) applies for the specific range of values $d \in [\frac{1}{2}, \frac{3}{2})$ and this can be extended by repeated use of partial summation in the definition. Model (1) directly provides a valid model for all values of $d$. Some interest in (1) has already been shown in the literature [e.g., Marinucci and Robinson (2000) and Robinson and Marinucci (2001)].

**3. Local Whittle estimation: consistency for $d \leq 1$ and inconsistency for $d > 1$.** Local Whittle (Gaussian semiparametric) estimation was developed by Künsch ([1987](#)) and Robinson ([1995](#)). Specifically, it starts with the following Gaussian objective function, defined in terms of the parameter $d$ and $G$:

$$(6) \qquad Q_m(G, d) = \frac{1}{m} \sum_{j=1}^{m} \left[ \log(G \lambda_j^{-2d}) + \frac{\lambda_j^{2d}}{G} I_x(\lambda_j) \right],$$

where $m$ is some integer less than $n$. The local Whittle procedure estimates $G$ and $d$ by minimizing $Q_m(G, d)$, so that

$$(\widehat{G}, \widehat{d}) = \underset{G \in (0, \infty), d \in [\Delta_1, \Delta_2]}{\arg \min} Q_m(G, d),$$

where $\Delta_1$ and $\Delta_2$ are numbers such that $-\frac{1}{2} < \Delta_1 < \Delta_2 < \infty$. It will be convenient in what follows to distinguish the true values of the parameters by the notation $G_0 = f_u(0)$ and $d_0$. Concentrating ([6](#)) with respect to $G$ as in Robinson ([1995](#)) gives

$$\widehat{d} = \underset{d \in [\Delta_1, \Delta_2]}{\arg \min} R(d),$$

where

$$R(d) = \log \widehat{G}(d) - 2d \frac{1}{m} \sum_{1}^{m} \log \lambda_j,$$

$$\widehat{G}(d) = \frac{1}{m} \sum_{1}^{m} \lambda_j^{2d} I_x(\lambda_j).$$

We now introduce the assumptions on $m$ and the stationary component $u_t$ in ([1](#)).

ASSUMPTION 1.

$$f_u(\lambda) \sim f_u(0) \in (0, \infty) \qquad \text{as } \lambda \to 0 + .$$

ASSUMPTION 2. In a neighborhood $(0, \delta)$ of the origin, $f_u(\lambda)$ is differentiable and

$$\frac{d}{d\lambda} \log f_u(\lambda) = O(\lambda^{-1}) \qquad \text{as } \lambda \to 0 + .$$

ASSUMPTION 3.

$$(7) \qquad u_t = C(L)\varepsilon_t = \sum_{j=0}^{\infty} c_j \varepsilon_{t-j}, \qquad \sum_{j=0}^{\infty} c_j^2 < \infty,$$

where $E(\varepsilon_t | F_{t-1}) = 0$, $E(\varepsilon_t^2 | F_{t-1}) = 1$ a.s., $t = 0, \pm 1, \ldots$, in which $F_t$ is the $\sigma$-field generated by $\varepsilon_s$, $s \le t$, and there exists a random variable $\varepsilon$ such that $E\varepsilon^2 < \infty$ and for all $\eta > 0$ and some $K > 0$, $\Pr(|\varepsilon_t| > \eta) \le K \Pr(|\varepsilon| > \eta)$.

ASSUMPTION 4.

$$\frac{1}{m} + \frac{m}{n} \to 0 \qquad \text{as } n \to \infty.$$

Assumptions 1–3 are analogous to Assumptions A1–A3 of Robinson (1995). However, we impose them in terms of $u_t$ rather than $X_t$. Assumption 4 is the same as Assumption A4 of Robinson (1995).

Lemma A.1(a) in Appendix A gives the following expression for $w_x(\lambda_s)$:

$$
\begin{aligned}
w_x(\lambda_s) = {} & \frac{D_n(e^{i\lambda_s}; \theta)}{1 - e^{i\lambda_s}} w_u(\lambda_s) \\
& - \frac{e^{i\lambda_s}}{1 - e^{i\lambda_s}} \frac{X_n - X_0}{\sqrt{2\pi n}} - \frac{1}{1 - e^{i\lambda_s}} \frac{\widetilde{U}_{\lambda_s n}(\theta)}{\sqrt{2\pi n}}.
\end{aligned}
\tag{8}
$$

Neglecting the third term of (8) as a remainder, $w_x(\lambda_s)$ is seen to comprise two terms—a function of the d.f.t. of $u_t$ and a function of $X_n$. As the value of $d$ changes, the stochastic magnitude of the two components changes, and this influences the asymptotic behavior of $w_x(\lambda_s)$. When $d < 1$, the first term dominates the second term and $w_x(\lambda_s)$ behaves like $\lambda_s^{-d} w_u(\lambda_s)$, being asymptotically uncorrelated for different frequencies. When $d > 1$, the second term becomes dominant and $w_x(\lambda_s)$ behaves like $\lambda_s^{-1}(X_n - X_0)/\sqrt{2\pi n}$, being perfectly correlated across all $\lambda_s$. This switching behavior of $w_x(\lambda_s)$ at $d = 1$ is a key determinant of the asymptotic properties of the local Whittle estimator, as well as other procedures like LP regression. When $d = 1$, the two terms have the same stochastic order and this leads to a form of asymptotic behavior that is particular to this case.

Theorem 3.1 below establishes that $\hat{d}$ is consistent for $d_0 \in (\frac{1}{2}, 1]$ and hence consistency carries over to the unit root case. While $\hat{G}$ is consistent for $d_0 \in (\frac{1}{2}, 1)$, however, it is inconsistent and tends to a random quantity when $d_0 = 1$.

THEOREM 3.1. *Suppose $X_t$ is generated by (1) with $d_0 \in [\Delta_1, \Delta_2]$ and Assumptions 1–4 hold. Then, for $d_0 \in (\frac{1}{2}, 1]$, $\hat{d} \overset{p}{\to} d_0$ as $n \to \infty$, and*

$$
\hat{G}(\hat{d}) \overset{d}{\to} \begin{cases} G_0, & \text{for } d_0 \in (\frac{1}{2}, 1), \\ G_0(1 + \chi_1^2), & \text{for } d_0 = 1. \end{cases}
$$

When $d_0 > 1$, $\hat{d}$ manifests very different behavior. It converges to unity in probability and the local Whittle estimator becomes inconsistent. So the local Whittle estimator is biased downward even in very large samples whenever the true value of $d$ is greater than unity. Kim and Phillips (1999) showed that the LP regression estimator also converges to unity when $d_0 > 1$.

THEOREM 3.2. *Under the same conditions as Theorem* 3.1, *for* $d_0 \in (1, M]$ *with* $1 < M < \infty$, $\widehat{d} \xrightarrow{p} 1$ *as* $n \to \infty$.

REMARK 3.3. Velasco (1999) showed that $\widehat{d}$ is consistent for $d_0 \in (\frac{1}{2}, 1)$ using the model (5). We conjecture that our consistency and inconsistency results for the local Whittle estimator for $d_0 = 1$ and $d_0 \in (1, \frac{3}{2})$ continue to hold under (5).

## 4. Local Whittle estimation: asymptotic distribution.
We introduce some further assumptions that are used in the results of this section.

ASSUMPTION 1′. For some $\beta \in (0, 2]$,

$$f_u(\lambda) = f_u(0)(1 + O(\lambda^\beta)), \qquad f_u(0) \in (0, \infty) \text{ as } \lambda \to 0+.$$

ASSUMPTION 2′. In a neighborhood $(0, \delta)$ of the origin, $C(e^{i\lambda})$ is differentiable and

$$\frac{d}{d\lambda} C(e^{i\lambda}) = O(\lambda^{-1}) \qquad \text{as } \lambda \to 0+.$$

ASSUMPTION 3′. Assumption 3 holds and also

$$E(\varepsilon_t^3 | F_{t-1}) = \mu_3,$$
$$E(\varepsilon_t^4 | F_{t-1}) = \mu_4 \qquad \text{a.s., } t = 0, \pm 1, \ldots,$$

for finite constants $\mu_3$ and $\mu_4$.

ASSUMPTION 4′. As $n \to \infty$,

$$\frac{1}{m} + \frac{m^{1+2\beta}(\log m)^2}{n^{2\beta}} \to 0.$$

ASSUMPTION 5′. Uniformly in $k = 0, 1, \ldots,$

$$\sum_{j \ge k} \gamma_j = O((\log(k+1))^{-4}), \qquad \sum_{j \ge k} c_j = O((\log(k+1))^{-4}),$$

$$\gamma_j \equiv E u_t u_{t+j}.$$

ASSUMPTION 6′. For the same $\beta \in (0, 2]$ as in Assumption 1′ and $\lambda, \lambda' \in (-\delta, \delta)$,

$$|C(e^{i\lambda}) - C(e^{i\lambda'})| \le C|\lambda - \lambda'|^{\min\{\beta, 1\}}, \qquad C \in (0, \infty).$$

Assumptions $1'$–$4'$ are analogous to Assumptions A1$'$–A4$'$ of Robinson ([1995]), except that our assumptions are in terms of $u_t$ rather than $X_t$. When $d_0 \in (\frac{1}{2}, 1)$, we need an additional assumption, Assumption $5'$, that controls the behavior of the tail sum of $c_j$ and $\gamma_j$. This assumption seems to be fairly mild. For instance, consider the stationary Gegenbauer process proposed by Gray, Zhang and Woodward ([1989]):

$$u_t = (1 - 2aL + L^2)^{-b} \varepsilon_t = C(L)\varepsilon_t,$$

$$t = 0, \pm 1, \pm 2, \ldots,$$

with $|a| < 1$ and $b \in (0, \frac{1}{2})$. Its spectral density is $f_u(\lambda) = \{4(\cos\lambda - a)^2\}^{-b}/2\pi$, which has a fractional pole at $\lambda_0 = \cos^{-1} a$. The asymptotic approximations for $c_j$ and $\gamma_j$ are given by [Gray, Zhang and Woodward ([1989]), pages 236–238]

(9)
$$c_j \sim \Lambda_1(a,b) \cos\{(j+b)\lambda_0 - b\pi/2\} j^{b-1},$$

$$\gamma_j \sim \Lambda_2(a,b) j^{2b-1} \sin(\pi b - j\lambda_0),$$

as $j \to \infty$, where $\Lambda_1(a,b)$ and $\Lambda_2(a,b)$ do not depend on $j$. Since $c_j$ and $\rho_j$ satisfy Assumption $5'$ [Zygmund ([1959]), Theorem 2.2, page 3], Assumption $5'$ allows for a pole and discontinuity in $f_u(\lambda)$ at $\lambda \neq 0$. However, Assumption $5'$ is not satisfied if $\gamma_k = (k+1)^{-1}(\log(k+1))^{-4}$. When $d_0 = 1$, Assumption $5'$ is not necessary, but instead we need Assumption $6'$. It requires $C(e^{i\lambda})$ to be Lip$(\min\{\beta, 1\})$ in the neighborhood of the origin.

The following theorems establish the asymptotic distribution of the local Whittle estimator for $d_0 \in (\frac{1}{2}, 1]$. When $d_0 \in (\frac{1}{2}, \frac{3}{4})$, $\hat{d}$ is asymptotically normally distributed, but $\hat{d}$ has a nonnormal limit distribution and slower rate of convergence when $d_0 \in [\frac{3}{4}, 1)$. This phenomenon occurs because, when $d_0$ is large, the stochastic magnitude of $X_n$ in the representation (8) becomes so large that it dominates the behavior of $\hat{d}$.

THEOREM 4.1. *Suppose $X_t$ is generated by* (1) *with $d_0 \in (\Delta_1, \Delta_2)$ and Assumptions $1'$–$5'$ hold. Then*

$$m^{1/2}(\hat{d} - d_0) \xrightarrow{d} \tfrac{1}{2}U \qquad \text{for } d_0 \in (\tfrac{1}{2}, \tfrac{3}{4}),$$

$$m^{1/2}(\hat{d} - d_0) \xrightarrow{d} \tfrac{1}{2}U + J(d_0)W^2 \qquad \text{for } d_0 = \tfrac{3}{4},$$

$$m^{2-2d_0}(\hat{d} - d_0) \xrightarrow{d} J(d_0)W^2 \qquad \text{for } d_0 \in (\tfrac{3}{4}, 1),$$

*where $J(d_0) = (2\pi)^{2d_0-2}\Gamma(d_0)^{-2}(2d_0-1)^{-3}(1-d_0)$ and $U$ and $W$ are mutually independent $N(0,1)$ random variables.*

When $d_0 = 1$ the two main components of $w_x(\lambda_s)$, that is, $w_u(\lambda_s)$ and $X_n/\sqrt{2\pi n}$, have the same stochastic magnitude, and the limit distribution of the local Whittle estimator turns out to be mixed normal (denoted as MN). Intriguingly, the variance of $\widehat{d}$ becomes smaller than the case where $d_0 < 1$, as was found in the corresponding case for LP regression [Phillips (1999b)].

THEOREM 4.2. *Suppose $X_t$ is generated by* (1) *with $d_0 = 1 \in (\Delta_1, \Delta_2)$ and Assumptions $1'$–$4'$ and $6'$ hold. Then*

$$m^{1/2}(\widehat{d} - d_0) \overset{d}{\to} \mathrm{MN}(0, \sigma^2(W))$$

$$\equiv \int_{-\infty}^{\infty} N(0, \sigma^2(h))\phi(h)\, dh,$$

*where $W$ is $N(0, 1)$, $\phi(\cdot)$ is standard normal p.d.f. and*

$$\sigma^2(h) = \frac{1}{4}\frac{1 + 2h^2}{1 + 2h^2 + h^4}.$$

REMARK 4.3. (a) When $d_0 = 1$, the variance of the limit distribution of $m^{1/2}(\widehat{d} - d_0)$ is less than $\frac{1}{4}$ since $\sigma^2(h) \leq \frac{1}{4}$ almost surely. Numerical evaluation gives

$$\sigma_d^2 = \frac{1}{4}\int_{-\infty}^{\infty} \frac{1 + 2h^2}{1 + 2h^2 + h^4}\frac{1}{\sqrt{2\pi}}\exp\left(-\frac{h^2}{2}\right)dh$$

$$= 0.2028.$$

Thus, the limit distribution of the local Whittle estimator has less dispersion when $d_0 = 1$ than it does in the stationary and $d_0 \in (\frac{1}{2}, \frac{3}{4})$ cases. A similar phenomenon applies in the limit theory for LP regression, where again the limit distribution is mixed normal when $d_0 = 1$ [Phillips (1999b)].

(b) Velasco (1999) shows asymptotic normality of the estimator for $d_0 \in (\frac{1}{2}, \frac{3}{4})$ using the model (5). We conjecture that the estimator has the same asymptotic distributions as those given above for $d_0 \in (\frac{3}{4}, 1]$ under (5), possibly with different $J(d_0)$, although the limit distribution for $d_0 = \frac{3}{4}$ might be difficult to derive.

**5. Fractional integration with a polynomial time trend.** In many applications, a nonstationary process is accompanied by a deterministic time trend. Accordingly, this section extends the analysis above to fractional processes with an $\alpha$-order ($\alpha > 0$) polynomial deterministic time trend. Specifically, the process $X_t$ is generated by the model

$$X_t = X_t^0 + X_0 + \mu t^\alpha = (1 - L)^{-d}u_t + X_0 + \mu t^\alpha$$

$$(10) \qquad = \sum_{k=0}^{t-1} \frac{(d)_k}{k!}u_{t-k} + X_0 + \mu t^\alpha, \qquad t = 0, 1, 2, \ldots, \ \mu \neq 0,$$

where $X_0$ and $u_t$ are defined as above. As shown in Appendix A, the d.f.t. of a time trend takes the following form, uniformly for $1 \le s \le m$ with $m = o(n)$:

$$w_{t^\alpha}(\lambda_s) = \frac{1}{\sqrt{2\pi n}} \sum_{t=1}^{n} t^\alpha e^{it\lambda_s}$$

$$= -\frac{1}{1 - e^{i\lambda_s}} \frac{n^\alpha}{\sqrt{2\pi n}}[1 + o(1)].$$

[See also Corbae, Ouliaris and Phillips (2002), who give exact formulae for d.f.t.'s of a time trend when $\alpha$ is a positive integer.] Therefore, neglecting the remainder term and $\widetilde{U}_{\lambda_s n}(\theta)$, we obtain the following expression of $w_x(\lambda_s)$:

$$
\begin{aligned}
(11) \quad w_x(\lambda_s) &\simeq -\frac{\mu}{1 - e^{i\lambda_s}} \frac{n^\alpha}{\sqrt{2\pi n}} + \frac{D_n(e^{i\lambda_s}; \theta)}{1 - e^{i\lambda_s}} w_u(\lambda_s) \\
&\quad - \frac{e^{i\lambda_s}}{1 - e^{i\lambda_s}} \frac{X_n^0 - X_0}{\sqrt{2\pi n}} \\
&\simeq C\mu\lambda_s^{-1} n^{\alpha - 1/2} \\
&\quad + O_p(\lambda_s^{-d}) + O_p(\lambda_s^{-1} n^{d-1}).
\end{aligned}
$$

When $\alpha > \frac{1}{2}$, the second term in (11) is dominated by either the first term (if $\alpha > d - \frac{1}{2}$) or the third term (if $\frac{1}{2} < \alpha < d - \frac{1}{2}$), and then $w_x(\lambda_s)$ behaves like $C(n)\lambda_s^{-1}$, where $C(n)$ does not depend on $s$. As a result, $\widehat{d}$ converges to unity in probability, and the local Whittle estimator is inconsistent except when the true value $d_0 = 1$. Since $X_n^0 = O_p(n^{d-1/2})$, this result might be regarded as an instance of a deterministic trend dominating a stochastic trend when $\alpha > d - \frac{1}{2}$. In the present case, because the d.f.t. of a deterministic trend is governed by the final observation, $n^\alpha$, the outcome for unfiltered, untapered data is the inconsistency of $\widehat{d}$. In consequence, some caution is needed in applying the Whittle estimator to investigate the degree of long range dependence when a time series exhibits trending behavior involving a deterministic trend of uncertain order. The same result holds if the deterministic trend $k_t$ is fractionally integrated in the sense that $(1 - L)^\alpha k_t = I\{t \ge 1\}$, because then $k_n \sim \Gamma(\alpha + 1)^{-1} n^\alpha$, as shown in Appendix A.

THEOREM 5.1. *Suppose $X_t$ is generated by (10) with $d_0 \in [\Delta_1, \Delta_2]$, $\alpha > \frac{1}{2}$, and Assumptions 1–4 hold. Then, for $d_0 \in (\frac{1}{2}, M]$ with $1 < M < \infty$, $\widehat{d} \xrightarrow{p} 1$ as $n \to \infty$.*

**6. Simulations and an empirical application.** First we report simulations that were conducted to examine the finite sample performance of the local Whittle estimator using (1) with $u_t \sim$ i.i.d. $N(0,1)$. All the results are based on 10,000 replications.

TABLE 1
*Simulation results for $d = 0.7$ and $d = 1.0$*

|  | $d = 0.7$ | | | $d = 1.0$ | | |
|---|---|---|---|---|---|---|
| $n$ | bias | s.d. | t.s.d.* | bias | s.d. | t.s.d. |
| 200 | 0.0002 | 0.1977 | 0.1336 | −0.0235 | 0.1779 | 0.1204 |
| 500 | 0.0093 | 0.1451 | 0.1066 | −0.0129 | 0.1280 | 0.0960 |
| 1000 | 0.0101 | 0.1162 | 0.0898 | −0.0102 | 0.1019 | 0.0809 |

*t.s.d. denotes theoretical standard deviation.

Table 1 shows the simulation results for $d = 0.7$ and $d = 1.0$. The sample size and $m$ were chosen to be $n = 200, 500, 1000$ and $m = [n^{0.5}]$. The estimator is seen to have smaller standard deviation when $d = 1.0$, corroborating the asymptotic theory.

Figure 1 plots the empirical distribution of the estimator for $d = 0.7$, $0.9$, $1.0$, $1.5$ when $n = 500$ and $m = [n^{0.5}]$. The estimator appears to have a symmetric distribution when $d \leq 1$, and the positive bias and skewness of the limit distribution for $d = 0.9$ are not evident for this sample size. When $d > 1$, distribution of the estimator is concentrated around unity, again corroborating the asymptotic result.

As an empirical illustration, the local Whittle estimator was applied to the historical economic time series considered in Nelson and Plosser (1982) and extended by Schotman and van Dijk (1991). We also estimate $d$ by first taking differences of the data, estimating $d − 1$ and adding unity to the estimate $\widehat{d − 1}$. This procedure is consistent for $\frac{1}{2} < d < 2$ and invariant to a linear trend. Table 2 shows the estimates based on both $m = n^{0.5}$ and

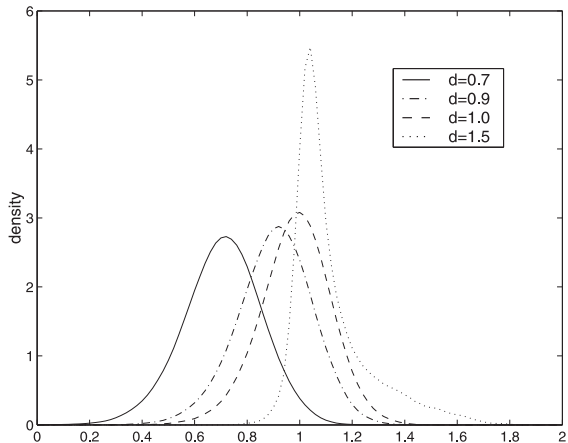

FIG. 1.  *Densities of the local Whittle estimator: $n = 500$, $m = n^{0.5}$.*

Table 2
*Estimates of d for US economic data*

| | $n$ | $m = n^{0.5}$ | | $m = n^{0.6}$ | |
|---|---|---|---|---|---|
| | | $d_{\mathrm{LW}}$ | $d_{\mathrm{LWD}}$ | $d_{\mathrm{LW}}$ | $d_{\mathrm{LWD}}$ |
| Real GNP | 62 | 0.990 | 0.626 | 0.946 | 0.719 |
| Nominal GNP | 62 | 0.983 | 0.901 | 0.930 | 0.909 |
| Real per capita GNP | 62 | 0.976 | 0.631 | 0.912 | 0.728 |
| Industrial production | 111 | 0.918 | 0.516 | 0.968 | 0.593 |
| Employment | 81 | 1.001 | 0.660 | 0.977 | 0.713 |
| Unemployment rate | 81 | 0.507 | 0.527 | 0.705 | 0.741 |
| GNP deflator | 82 | 1.143 | 0.973 | 1.049 | 1.099 |
| CPI | 111 | 1.020 | 1.227 | 0.828 | 1.176 |
| Nominal wage | 71 | 1.080 | 1.026 | 1.015 | 0.983 |
| Real wage | 71 | 1.105 | 0.785 | 1.030 | 0.822 |
| Money stock | 82 | 1.042 | 0.913 | 0.993 | 1.232 |
| Velocity of money | 102 | 1.055 | 0.932 | 0.970 | 0.782 |
| Bond yield | 71 | 0.676 | 1.261 | 0.740 | 1.370 |
| Stock prices | 100 | 0.914 | 0.860 | 0.984 | 0.755 |

$m = n^{0.6}$. These series produce long memory estimates over a wide interval that ranges from around 0.5 for the unemployment rate to 1.38 for the bond yield. For the unemployment rate, the local Whittle estimate from the raw data ($\hat{d}_{\mathrm{LW}}$) and the local Whittle estimate from the differenced data ($\hat{d}_{\mathrm{LWD}}$) are very close together, both indicating only marginal nonstationarity in the data. For the bond yield, $\hat{d}_{\mathrm{LWD}}$ is very different from $\hat{d}_{\mathrm{LW}}$. Especially for the GNP measures, industrial production and employment, the presence of a linear trend component in the data [which is supported by much of the empirical work with this data set following Nelson and Plosser (1982)] appears to bias $\hat{d}_{\mathrm{LW}}$ heavily toward unity. These particular results indicate that, although the local Whittle estimator is consistent for $\frac{1}{2} < d \leq 1$, the use of differenced data or even data tapering [Velasco (1999) and Hurvich and Chen (2000)] may be preferable, unless the time series clearly does not involve a deterministic trend and values of $d > 1$ are not suspected.

**7. Concluding remarks.** The results of the present paper have a negative character, revealing that the local Whittle estimator is not a good general purpose estimator when the value of $d$ may take on values in the nonstationary zone beyond $\frac{3}{4}$. The asymptotic theory is discontinuous at $d = \frac{3}{4}$ and again at $d = 1$, it is awkward to use and the estimator is inconsistent beyond unity.

This paper has not explicitly addressed the issue of what semiparametric estimation procedure is a good general purpose procedure for possibly nonstationary cases. Data differencing and data tapering have been explored [Velasco (1999) and Hurvich and Chen (2000)], are easy to implement and have been shown to extend the range of applicability of the local Whittle estimator. However, these approaches do have some disadvantages, such as the need to determine the appropriate order of differencing and the effects of tapering on variance. Another approach is to use the exact form of the local Whittle estimator suggested in Phillips (1999a), which does not rely on differencing or tapering. This estimator has recently been shown by the authors [Shimotsu and Phillips (2002)] to be consistent and to have the same $N(0, \frac{1}{4})$ limit distribution for all values of $d$. While it is still too early for a definitive answer to the question of what is a good general purpose semi-parametric estimator of $d$ that allows for nonstationarity, these approaches offer some useful alternatives for applied researchers, and the present paper is at least a cautionary tale about performance characteristics of the local Whittle estimator in the nonstationary environment.

## APPENDIX A

**Technical lemmas.** In this and the following sections, $x^*$ denotes the complex conjugate of $x$, and $|x|_+$ denotes $\max\{x, 1\}$.

LEMMA A.1 [Phillips (1999a), Theorems 2.2 and 2.7]. (a) *If $X_t$ follows* (1), *then*

$$
\begin{aligned}
(12) \quad w_x(\lambda)(1 - e^{i\lambda}) = {} & D_n(e^{i\lambda}; \theta) w_u(\lambda) \\
& - \frac{e^{in\lambda}}{\sqrt{2\pi n}} \widetilde{U}_{\lambda n}(\theta) - \frac{e^{i\lambda}}{\sqrt{2\pi n}} (e^{in\lambda} X_n - X_0),
\end{aligned}
$$

*where* $D_n(e^{i\lambda}; \theta) = \sum_{k=0}^{n} \frac{(-\theta)_k}{k!} e^{ik\lambda}$, $\theta = 1 - d$ *and*

$$
\begin{aligned}
(13) \quad & \widetilde{U}_{\lambda n}(\theta) = \widetilde{D}_{n\lambda}(e^{-i\lambda} L; \theta) u_n = \sum_{p=0}^{n-1} \widetilde{\theta}_{\lambda p} e^{-ip\lambda} u_{n-p}, \\
& \widetilde{\theta}_{\lambda p} = \sum_{k=p+1}^{n} \frac{(-\theta)_k}{k!} e^{ik\lambda}.
\end{aligned}
$$

(b) *If $X_t$ follows* (1) *with $d = 1$, then*

$$
(14) \quad w_x(\lambda)(1 - e^{i\lambda}) = w_u(\lambda) - \frac{e^{i\lambda}}{\sqrt{2\pi n}} (e^{in\lambda} X_n - X_0).
$$

LEMMA A.2. *For $\theta > -1$ and uniformly in $s = 1, 2, \ldots, m$ with $m = o(n)$,*

$$
(15) \quad D_n(e^{i\lambda_s}; \theta) = (1 - e^{i\lambda_s})^\theta + O(n^{-\theta} s^{-1}).
$$

PROOF. We have

$$
D_n(e^{i\lambda_s};\theta) = \sum_0^\infty \frac{(-\theta)_k}{k!} e^{ik\lambda_s} - \sum_{n+1}^\infty \frac{(-\theta)_k}{k!} e^{ik\lambda_s}
$$

(16)

$$
= {}_2F_1(-\theta,1;1;e^{i\lambda_s}) - \sum_{n+1}^\infty \frac{k^{-\theta-1}}{\Gamma(-\theta)} e^{ik\lambda_s} + O\left(\sum_{n+1}^\infty k^{-\theta-2}\right),
$$

since $(-\theta)_k/k! = \Gamma(-\theta)^{-1} k^{-\theta-1}(1 + O(k^{-1}))$ [Erdélyi [1953], page 47]. Because $\theta > -1$ and $s \neq 0$, the first term in (16) converges and equals to $(1 - e^{i\lambda_s})^\theta$ [Erdélyi [1953], page 57]. For the second term in (16), by Theorem 2.2 of Zygmund [[1959], page 3] we have

$$
\left| \sum_{n+1}^\infty k^{-\theta-1} e^{ik\lambda_s} \right|
$$

$$
\leq (n+1)^{-\theta-1} \max_N \left| \sum_{n+1}^{n+N} e^{ik\lambda_s} \right|
$$

$$
= O(n^{-\theta} s^{-1}).
$$

The third term in (16) is necessarily $O(n^{-\theta} s^{-1})$ because $\sum_{n+1}^\infty k^{-\theta-2} = O(n^{-\theta-1})$. □

LEMMA A.3. (a) $\lambda^{-\theta}(1 - e^{i\lambda})^\theta = e^{-(\pi/2)\theta i} + O(\lambda)$ as $\lambda \to 0+$.
(b) For $\theta > -1$ and uniformly in $s = 1, 2, \ldots, m$ with $m = o(n)$,

(17) $$\lambda_s^{-\theta} D_n(e^{i\lambda_s};\theta) = e^{-(\pi/2)\theta i} + O(\lambda_s) + O(s^{-1-\theta}).$$

PROOF. For (a), since $|1 - e^{i\lambda}| = |2\sin(\lambda/2)|$ and $\arg(1 - e^{i\lambda}) = (\lambda - \pi)/2$ for $0 \leq \lambda < \pi$, we can write $(1 - e^{i\lambda})^\theta$ in polar form as $|2\sin(\lambda/2)|^\theta \exp[i\theta(\lambda - \pi)/2]$. It follows that

$$
\lambda^{-\theta}(1 - e^{i\lambda})^\theta = \lambda^{-\theta}(\lambda + O(\lambda^3))^\theta [\exp(-i\theta\pi/2) + O(\lambda)]
$$

$$
= e^{-(\pi/2)\theta i} + O(\lambda)
$$

giving the stated result. Statement (b) follows from (a) and Lemma A.2. □

LEMMA A.4. Uniformly in $p = 0, 1, \ldots, n-1$ and $s = 1, 2, \ldots, m$ with $m = o(n)$:

(a)

(18) $$\widetilde{\theta}_{\lambda_s p} = \begin{cases} O(|p|_+^{-\theta}) = O(|p|_+^{d-1}), & \text{for } \theta > 0, \\ O(n^{-\theta}) = O(n^{d-1}), & \text{for } \theta \in (-1, 0), \end{cases}$$

(b)

$$\widetilde{\theta}_{\lambda_s p} = O(|p|_+^{-\theta-1} n s^{-1})$$

(19)

$$= O(|p|_+^{d-2} n s^{-1}), \qquad \text{for } \theta > -1.$$

PROOF.   Observe

$$\widetilde{\theta}_{\lambda_s p} = \Gamma(-\theta)^{-1} \sum_{p+1}^{n} k^{-\theta-1} e^{ik\lambda_s} + O\left(\sum_{p+1}^{n} k^{-\theta-2}\right).$$

The required results follow from

$$\sum_{p+1}^{n} k^{-\theta-1} = \begin{cases} O(|p|_+^{-\theta}), & \text{for } \theta > 0, \\ O(n^{-\theta}), & \text{for } \theta \in (-1, 0), \end{cases}$$

$$\sum_{p+1}^{n} k^{-\theta-2} \le \sum_{p+1}^{n} k^{-\theta-1},$$

$$\left| \sum_{p+1}^{n} k^{-\theta-1} e^{ik\lambda_s} \right| \le (p+1)^{-\theta-1} \max_{N} \left| \sum_{p+1}^{p+N} e^{ik\lambda_s} \right|$$

$$= O(|p|_+^{-\theta-1} n s^{-1})$$

and $\sum_{p+1}^{n} k^{-\theta-2} = O(|p|_+^{-\theta-1})$.   □

LEMMA A.5.   (a) *Under the assumptions of Theorem* 3.1, *we have the following:*

(a1)                $E|\widetilde{U}_{\lambda_s n}(\theta)|^2 = O(h_{ns}(\theta))$,

(a2)     $E(X_n - X_0 - C(1) X_n^\varepsilon)^2 = o(n^{2d-1})$,

*uniformly in* $s = 1, 2, \ldots, m$, *where* $X_n^\varepsilon = \sum_{k=0}^{n-1} \frac{(d)_k}{k!} \varepsilon_{n-k}$ *and*

$$h_{ns}(\theta) = \begin{cases} n^{1-2\theta} s^{2\theta-1} = n^{2d-1} s^{1-2d}, & \text{for } \theta \in (-\frac{1}{2}, \frac{1}{2}), \\ n^{1-2\theta} s^{2\theta-1} (\log(s+1))^2 & \\ \quad = n^{2d-1} s^{1-2d} (\log(s+1))^2, & \text{for } \theta = -\frac{1}{2}. \end{cases}$$

(b) *Under the assumptions of Theorem* 4.1, *we have, uniformly in* $s = 1, \ldots, m$,

$$E|\widetilde{U}_{\lambda_s n}(\theta) - C(1)\widetilde{\varepsilon}_{\lambda_s n}(\theta)|^2$$
$$= O(n^{1-2\theta} s^{2\theta-1} (\log n)^{-4} + n^{1-2\theta} s^{-2}) \qquad \text{for } \theta \in (-\tfrac{1}{2}, \tfrac{1}{2}).$$

PROOF. (a) We prove (a1) first. When $\theta = 0$ the stated result follows because $\widetilde{U}_{\lambda_s n}(\theta) = 0$. When $\theta \neq 0$ define $a_p = \bar{\theta}_{\lambda_s p} e^{-ip\lambda_s}$ so that $\widetilde{U}_{\lambda_s n}(\theta) = \sum_{p=0}^{n-1} a_p u_{n-p}$. We suppress the dependence of $a_p$ on $\theta$ and $\lambda_s$. Summation by parts gives

$$(20) \qquad \widetilde{U}_{\lambda_s n}(\theta) = \sum_{p=0}^{n-2}(a_p - a_{p+1})\sum_{j=0}^{p} u_{n-j} + a_{n-1}\sum_{j=0}^{n-1} u_{n-j}.$$

Observe that

$$a_p - a_{p+1}$$

$$= \sum_{k=p+1}^{n} \frac{(-\theta)_k}{k!} e^{i(k-p)\lambda_s} - \sum_{k=p+2}^{n} \frac{(-\theta)_k}{k!} e^{i(k-p-1)\lambda_s}$$

$$= \sum_{k=p+1}^{n} \frac{(-\theta)_k}{k!} e^{i(k-p)\lambda_s} - \sum_{l=p+1}^{n-1} \frac{(-\theta)_{l+1}}{(l+1)!} e^{i(l-p)\lambda_s}$$

$$= \sum_{k=p+1}^{n-1} \left[\frac{(-\theta)_k}{k!} - \frac{(-\theta)_{k+1}}{(k+1)!}\right] e^{i(k-p)\lambda_s} + \frac{(-\theta)_n}{n!} e^{-ip\lambda_s}$$

$$= \sum_{k=p+1}^{n-1} \frac{(1+\theta)\Gamma(k-\theta)}{\Gamma(-\theta)\Gamma(k+2)} e^{i(k-p)\lambda_s} + \frac{(-\theta)_n}{n!} e^{-ip\lambda_s},$$

where the fourth line follows from

$$\frac{(-\theta)_k}{k!} - \frac{(-\theta)_{k+1}}{(k+1)!}$$

$$= -\frac{\Gamma(k-\theta)}{\Gamma(-\theta-1)\Gamma(k+2)} = \frac{(1+\theta)\Gamma(k-\theta)}{\Gamma(-\theta)\Gamma(k+2)}.$$

Define

$$b_{np} = \sum_{k=p+1}^{n-1} \frac{(1+\theta)\Gamma(k-\theta)}{\Gamma(-\theta)\Gamma(k+2)} e^{i(k-p)\lambda_s},$$

and then, since $a_{n-1} = ((-\theta)_n/n!)e^{-i(n-1)\lambda_s}$,

$$\widetilde{U}_{\lambda_s n}(\theta) = \sum_{p=0}^{n-2} b_{np} \sum_{j=0}^{p} u_{n-j} + \frac{(-\theta)_n}{n!} \sum_{p=0}^{n-2} e^{-ip\lambda_s} \sum_{j=0}^{p} u_{n-j}$$

$$+ \frac{(-\theta)_n}{n!} e^{-i(n-1)\lambda_s} \sum_{j=0}^{n-1} u_{n-j}$$

$$= \sum_{p=0}^{n-2} b_{np} \sum_{j=0}^{p} u_{n-j} + \frac{(-\theta)_n}{n!} \sum_{p=0}^{n-1} e^{-ip\lambda_s} \sum_{j=0}^{p} u_{n-j}$$

$$= U_{1n} + U_{2n}.$$

We proceed to show that the $U_{\cdot n}$ are of the stated order. First, for $U_{1n}$ we have

$$(21) \qquad b_{np} = O(\min\{|p|_+^{-\theta-1}, |p|_+^{-\theta-2} n s^{-1}\})$$

uniformly in $p = 0, \ldots, n-1$ and $s = 1, \ldots, m$. Equation (21) holds because

$$b_{np} = \frac{1+\theta}{\Gamma(-\theta)} e^{-ip\lambda_s} \sum_{k=p+1}^{n-1} k^{-\theta-2} e^{ik\lambda_s} + O\left( \sum_{p+1}^{n-1} k^{-\theta-3} \right)$$

and

$$\sum_{p+1}^{n-1} k^{-\theta-2} = O(|p|_+^{-\theta-1}),$$

$$\sum_{p+1}^{n-1} k^{-\theta-3} = O(|p|_+^{-\theta-2}),$$

$$\left| \sum_{p+1}^{n-1} k^{-\theta-2} e^{ik\lambda_s} \right| \leq (p+1)^{-\theta-2} \max_N \left| \sum_{p+1}^{p+N} e^{ik\lambda_s} \right|$$

$$= O(|p|_+^{-\theta-2} n s^{-1}).$$

Next,

$$(22) \qquad E\left( \sum_0^p u_{n-j} \right)^2$$

$$= (p+1) \sum_{-p}^{p} (1 - |j|/(p+1)) \gamma_j = O(|p|_+), \qquad \gamma_j \equiv E u_t u_{t+j},$$

for $p = 0, \ldots, n-1$, and it follows from Minkowski's inequality that

$$E|U_{1n}|^2 = O\left( \left( \sum_{p=0}^{n-2} |b_{np}| |p|_+^{1/2} \right)^2 \right)$$

$$= O\left( \left( \sum_{p=0}^{n/s} |p|_+^{-\theta-1/2} + \sum_{p=n/s}^{n} p^{-\theta-3/2} n s^{-1} \right)^2 \right)$$

$$= \begin{cases} O(n^{1-2\theta} s^{2\theta-1}), & \theta \in (-\frac{1}{2}, \frac{1}{2}), \\ O(n^{1-2\theta} s^{2\theta-1} (\log(s+1))^2), & \theta = -\frac{1}{2}. \end{cases}$$

For $U_{2n}$, we rewrite the sum as

$$U_{2n} = \frac{(-\theta)_n}{n!} \sum_{p=0}^{n-1} e^{-ip\lambda_s} \sum_{j=0}^{p} u_{n-j}$$

$$= \frac{(-\theta)_n}{n!} \sum_{n-p=1}^{n} e^{i(n-p)\lambda_s} \sum_{n-j=n-p}^{n} u_{n-j}$$

$$= \frac{(-\theta)_n}{n!} \sum_{k=1}^{n} u_k \sum_{q=1}^{k} e^{iq\lambda_s}$$

$$= \frac{(-\theta)_n}{n!} \sum_{k=1}^{n} u_k \frac{e^{i\lambda_s}(1 - e^{ik\lambda_s})}{1 - e^{i\lambda_s}}$$

$$= \frac{(-\theta)_n}{n!} \frac{e^{i\lambda_s}}{1 - e^{i\lambda_s}} \sum_{k=1}^{n} u_k - \frac{(-\theta)_n}{n!} \frac{e^{i\lambda_s}}{1 - e^{i\lambda_s}} (2\pi n)^{1/2} w_u(\lambda_s);$$

$E|U_{2n}|^2 = O(n^{1-2\theta}s^{-2})$ follows from (22) and $E|w_u(\lambda_s)|^2 = O(1)$ [Robinson (1995), page 1637], and the stated result follows because $s^{-2} \leq s^{2\theta-1}$.

We move to the proof of (a2). Define $a_p = (d)_p/p!$ so that $X_n = \sum_{p=0}^{n-1} a_p u_{n-p} + X_0$. Similar to the above, summation by parts gives

$$X_n - X_0 - C(1)X_n^\varepsilon$$

$$= \sum_{p=0}^{n-2} (a_p - a_{p+1}) \sum_{j=0}^{p} (u_{n-j} - C(1)\varepsilon_{n-j}) + a_{n-1} \sum_{j=0}^{n-1} (u_{n-j} - C(1)\varepsilon_{n-j}).$$

Since $a_p - a_{p+1} = -\frac{\Gamma(d+p)}{\Gamma(d-1)\Gamma(p+2)} = O(|p|_+^{d-2})$ and $a_p = O(|p|_+^{d-1})$, the stated result follows if

$$(23) \qquad E\left[\sum_{j=0}^{p} (u_{n-j} - C(1)\varepsilon_{n-j})\right]^2 = o(p) \qquad \text{as } p \to \infty.$$

Now

$$E\left[\sum_{j=0}^{p} u_{n-j}\right]^2 = \sum_{j=-p}^{p} (p + 1 - |j|)\gamma_j,$$

$$E\left[C(1) \sum_{j=0}^{p} \varepsilon_{n-j}\right]^2 = (p+1)C(1)^2 = (p+1) \sum_{j=-\infty}^{\infty} \gamma_j,$$

$$E\left[C(1) \sum_{j=0}^{p} u_{n-j} \sum_{l=0}^{p} \varepsilon_{n-l}\right] = C(1)E\left[\sum_{j=0}^{p} \sum_{r=0}^{\infty} c_r \varepsilon_{n-j-r} \sum_{l=0}^{p} \varepsilon_{n-l}\right]$$

$$= C(1) \sum_{j=0}^{p} \sum_{l=0}^{p} \sum_{r=0}^{\infty} c_r I\{r = l - j\}$$

$$= C(1) \sum_{r=0}^{p} (p + 1 - r) c_r,$$

and it follows that $E[\sum_{j=0}^{p}(u_{n-j} - C(1)\varepsilon_{n-j})]^2$ is equal to

$$
\begin{aligned}
(24) \qquad & -(p+1) \sum_{|j| \geq p+1} \gamma_j - 2 \sum_1^p j \gamma_j \\
& + 2C(1)(p+1) \sum_{r \geq p+1} c_r - 2C(1) \sum_1^p r c_r,
\end{aligned}
$$

which is $o(p)$ from $\sum_{-\infty}^{\infty} \gamma_j$, $\sum_0^{\infty} c_r < \infty$ and Kronecker's lemma. Therefore, (23) and the stated result follow.

(b) Let $M$ be a generic finite positive constant. We collect some facts that are used repeatedly: for $\alpha \in (-1, C)$ and $q \geq 2$,

$$(25) \qquad \sum_{l=2}^{q} (\log l)^{-4} \leq (\log 2)^{-4} \sum_2^{\sqrt{q}} + (\tfrac{1}{2} \log q)^{-4} \sum_{\sqrt{q}}^{q} \leq Mq(\log q)^{-4},$$

$$
\begin{aligned}
(26) \qquad & \sum_{l=0}^{q} |l|_+^{\alpha}(\log(l+2))^{-2} \leq (\log 2)^{-2} \sum_0^{\sqrt{q}} |l|_+^{\alpha} + (\tfrac{1}{2} \log q)^{-2} \sum_{\sqrt{q}}^{q} l^{\alpha} \\
& \leq Mq^{\alpha+1}(\log q)^{-2}.
\end{aligned}
$$

Proceeding similarly to the proof of (a1), we obtain

$$\widetilde{U}_{\lambda_s n}(\theta) - C(1)\widetilde{\varepsilon}_{\lambda_s n}(\theta) = \dot{U}_{1n} + \dot{U}_{2n},$$

where

$$
\begin{aligned}
\dot{U}_{1n} &= \sum_{p=0}^{n-2} b_{np} \sum_{j=0}^{p} (u_{n-j} - C(1)\varepsilon_{n-j}), \\
\dot{U}_{2n} &= \frac{(-\theta)_n}{n!} \frac{e^{i\lambda_s}}{1 - e^{i\lambda_s}} \sum_{k=1}^{n} (u_k - C(1)\varepsilon_k) \\
& \qquad - \frac{(-\theta)_n}{n!} \frac{e^{i\lambda_s}}{1 - e^{i\lambda_s}} (2\pi n)^{1/2} [w_u(\lambda_s) - C(1)w_{\varepsilon}(\lambda_s)]
\end{aligned}
$$

and $b_{np}$ is defined in (21).

First, we show that, uniformly in $p = 0, \ldots, n-1$,

$$(27) \quad E\left[\sum_{j=0}^{p}(u_{n-j} - C(1)\varepsilon_{n-j})\right]^2$$
$$= O(|p|_+(\log(p+2))^{-4}).$$

When $p = 0$, (28) follows immediately. When $p \geq 1$, from (24) the left-hand side of (27) is equal to

$$-(p+1)\sum_{|j|\geq p+1}\gamma_j - 2\sum_1^p j\gamma_j$$
$$+ 2C(1)(p+1)\sum_{r\geq p+1}c_r - 2C(1)\sum_1^p r c_r.$$

The first and third terms are bounded uniformly in $p$ by $p(\log(p+2))^{-4}$ from Assumption 5'. For the second term we have

$$\left|\sum_1^p j\gamma_j\right| = \left|\sum_{j=1}^p\sum_{k=j}^p\gamma_k\right|$$
$$= O\left(\sum_1^p(\log(j+1))^{-4}\right)$$
$$= O((p+1)(\log(p+1))^{-4})$$

uniformly in $p$, where the third equality follows from (25); $\sum_1^p r c_r = O((p+1) \times (\log(p+1))^{-4})$ follows from the same argument, and (27) follows.

From Minkowski's inequality, (21) and (27), $(E|\dot{U}_{1n}|^2)^{1/2}$ is bounded by

$$\sum_{p=0}^{n-2}|b_{np}||p|_+^{1/2}(\log(p+2))^{-2}$$
$$= O\left(\sum_{p=0}^{n/s}|p|_+^{-\theta-1/2}(\log(p+2))^{-2} + \sum_{p=n/s}^n|p|_+^{-\theta-3/2}ns^{-1}(\log(p+2))^{-2}\right)$$
$$= O\left((n/s)^{1/2-\theta}(\log(n/s))^{-2} + (\log(n/s))^{-2}ns^{-1}\sum_{n/s}^n|p|_+^{-\theta-3/2}\right)$$
$$= O(n^{1/2-\theta}s^{\theta-1/2}(\log n)^{-2}),$$

where the third line follows from (26), and the fourth line follows because $(\log(n/s))^{-2} \leq (\log(n/m))^{-2} = O((\log n)^{-2})$; $E|\dot{U}_{2n}|^2 = O(n^{1-2\theta}s^{-2})$ follows from (27) and $E|w_u(\lambda_s) - C(1)w_\varepsilon(\lambda_s)|^2 = O(1)$ [Robinson (1995), page 1637], giving the stated result. $\square$

LEMMA A.6. *Under the assumptions of Theorem 4.2, we have, for $j = 1, \ldots, m$,*

$$E|w_u(\lambda_j) - C(e^{i\lambda_j})w_\varepsilon(\lambda_j)|^2 = \begin{cases} O(n^{-\beta}), & \text{for } \beta \in (0,1), \\ O(n^{-1}\log n), & \text{for } \beta \in [1,2]. \end{cases}$$

PROOF. The proof essentially follows from Theorem 3.15 of Zygmund [(1959), page 91]. An elementary calculation gives

$$(28) \quad \begin{aligned} & Ew_u(\lambda_j)w_\varepsilon^*(\lambda_j) - C(e^{i\lambda_j})/2\pi \\ & = \frac{1}{2\pi} \int_{-\pi}^{\pi} [C(e^{i\lambda}) - C(e^{i\lambda_j})]K(\lambda - \lambda_j)\, d\lambda, \end{aligned}$$

where $K(\lambda) = (2\pi n)^{-1}\sum_1^n\sum_1^n e^{i(t-s)\lambda}$ is Fejér's kernel. From Zygmund [(1959), page 90], $|K(\lambda)| \leq An^{-1}\lambda^{-2}$ and $|K(\lambda)| \leq An$ for a finite constant $A$. Assumption 6' implies $|C(e^{i\lambda}) - C(e^{i\lambda_j})| \leq C|\lambda - \lambda_j|^{\min\{\beta,1\}}$ for $|\lambda - \lambda_j| \leq \delta/2$ and large enough $n$. Therefore, if we split the integral (28), each part is bounded as follows:

$$\begin{aligned} \int_{-\pi}^{\lambda_j-\delta/2} + \int_{\lambda_j+\delta/2}^{\pi} &= O\left(n^{-1}\int_{\delta/2}^{\pi}\lambda^{-2}\, d\lambda\right) \\ &= O(n^{-1}), \\ \int_{\lambda_j-\delta/2}^{\lambda_j-1/n} + \int_{\lambda_j+1/n}^{\lambda_j+\delta/2} &= O\left(n^{-1}\int_{1/n}^{\delta/2}\lambda^{\min\{\beta-2,-1\}}\, d\lambda\right) \\ &= \begin{cases} O(n^{-\beta}), & \text{for } \beta \in (0,1), \\ O(n^{-1}\log n), & \text{for } \beta \in (1,2] \end{cases} \end{aligned}$$

and

$$\begin{aligned} \int_{\lambda_j-1/n}^{\lambda_j+1/n} &= O\left(n\int_0^{1/n}\lambda^{\min\{\beta,1\}}\, d\lambda\right) \\ &= O(n^{-\min\{\beta,1\}}). \end{aligned}$$

Hence, $Ew_u(\lambda_j)w_\varepsilon^*(\lambda_j) - C(e^{i\lambda_j})/2\pi$ has the stated order; $EI_u(\lambda_j) - f_u(\lambda_j)$ has the same order by a similar argument, and the order of

$$\begin{aligned} & E|w_u(\lambda_j) - C(e^{i\lambda_j})w_\varepsilon(\lambda_j)|^2 \\ & = E[I_u(\lambda_j) - 2\operatorname{Re}[w_u(\lambda_j)C^*(e^{i\lambda_j})w_\varepsilon(\lambda_j)] + 2\pi f_u(\lambda_j)I_\varepsilon(\lambda_j)] \end{aligned}$$

follows.   □

LEMMA A.7. *Let $v_t = I\{t \geq 1\}$ and $\Delta^{-\alpha}v_t = (1-L)^{-\alpha}v_t$ with $\alpha > 0$.*
*Then uniformly in $1 \leq s \leq m$ with $m = o(n)$ the following hold:*

(a) $\quad w_{\Delta^{-\alpha}v}(\lambda_s) = -\dfrac{e^{i\lambda_s}}{1-e^{i\lambda_s}}\dfrac{1}{\Gamma(\alpha+1)}\dfrac{n^\alpha}{\sqrt{2\pi n}}[1 + O(s^{-\min\{\alpha,1\}})];$

(b) $\quad w_{t^\alpha}(\lambda_s) = -\dfrac{1}{1-e^{i\lambda_s}}\dfrac{n^\alpha}{\sqrt{2\pi n}}[1 + O(s^{-\min\{\alpha,1\}}) + O(n^{-1}s)].$

PROOF. For part (a), first consider the case $\alpha \in (0,1]$. From Lemma
A.1(b)

$$(29) \qquad \begin{aligned} w_{\Delta^{-\alpha}v}(\lambda_s) &= (1-e^{i\lambda_s})^{-1}w_{\Delta^{-\alpha+1}v}(\lambda_s) \\ &\quad - (1-e^{i\lambda_s})^{-1}e^{i\lambda_s}\Delta^{-\alpha}v_n/\sqrt{2\pi n}. \end{aligned}$$

For $\alpha = 1$, since $w_v(\lambda_s) = 0$ it follows that

$$\begin{aligned} w_{\Delta^{-1}v}(\lambda_s) &= -(1-e^{i\lambda_s})^{-1}e^{i\lambda_s}\Delta^{-1}v_n/\sqrt{2\pi n} \\ &= -(1-e^{i\lambda_s})^{-1}e^{i\lambda_s}n/\sqrt{2\pi n}. \end{aligned}$$

From

$$\frac{(d)_k}{k!} - \frac{(d)_{k-1}}{(k-1)!} = \frac{\Gamma(k-1+d)}{\Gamma(d-1)\Gamma(k+1)} = \frac{(d-1)_k}{k!}$$

and the fact that $(\alpha-1)_0/0! = (\alpha)_0/0! = 1$ we obtain

$$\begin{aligned} (1-L)^{-\alpha+1}v_t &= \sum_{k=0}^{t-1}\frac{(\alpha-1)_k}{k!} \\ &= \sum_{k=1}^{t-1}\left[\frac{(\alpha)_k}{k!} - \frac{(\alpha)_{k-1}}{(k-1)!}\right] + \frac{(\alpha)_0}{0!} = \frac{(\alpha)_{t-1}}{(t-1)!}. \end{aligned}$$

Hence, for $\alpha \in (0,1)$ from Lemma A.2 we have

$$\begin{aligned} w_{\Delta^{-\alpha+1}v}(\lambda_s) &= \frac{1}{\sqrt{2\pi n}}\sum_{t=1}^{n}\frac{(\alpha)_{t-1}}{(t-1)!}e^{it\lambda_s} \\ &= \frac{e^{i\lambda_s}}{\sqrt{2\pi n}}\left[D_n(e^{i\lambda_s};-\alpha) - \frac{(\alpha)_n}{n!}\right] \\ &= \frac{e^{i\lambda_s}}{\sqrt{2\pi n}}[(1-e^{i\lambda_s})^{-\alpha} + O(n^\alpha s^{-1})]. \end{aligned}$$

Then the stated result follows because

$$\Delta^{-\alpha}v_n = \frac{(\alpha+1)_{n-1}}{(n-1)!} = \frac{1}{\Gamma(\alpha+1)}n^\alpha[1 + O(n^{-1})],$$

so that the second term on the right-hand side of (29) dominates the first term. The result for $\alpha > 1$ is derived from (29) and by induction.

For part (b), observe that

$$
\begin{aligned}
w_{\Delta^{-\alpha}v}(\lambda_s) &= \frac{e^{i\lambda_s}}{\sqrt{2\pi n}} \sum_{t=0}^{n-1} \frac{(\alpha+1)_t}{t!} e^{it\lambda_s} \\
&= \frac{e^{i\lambda_s}}{\sqrt{2\pi n}} + \frac{e^{i\lambda_s}}{\sqrt{2\pi n}} \sum_1^{n-1} \left[ \frac{1}{\Gamma(\alpha+1)} t^\alpha + O(t^{\alpha-1}) \right] e^{it\lambda_s} \\
&= e^{i\lambda_s} \Gamma(\alpha+1)^{-1} w_{t^\alpha}(\lambda_s) + O(n^{\alpha-1/2}),
\end{aligned}
$$

and the required result follows from part (a).  $\square$

## APPENDIX B

**Proofs of theorems.**

PROOF OF THEOREM 3.1.   For notational simplicity we assume $X_0 = 0$ throughout the proof, but the result carries over for general $X_0$ with $X_n - X_0$ replacing $X_n$. We follow the approach developed by Robinson (1995) for the stationary case. Define $G(d) = G_0 m^{-1} \sum_1^m \lambda_j^{2d-2d_0}$ and $S(d) = R(d) - R(d_0)$. For arbitrarily small $\Delta > 0$, define $\Theta_1 = \{d : d_0 - \frac{1}{2} + \Delta \le d \le \Delta_2\}$ and $\Theta_2 = \{d : \Delta_1 \le d < d_0 - \frac{1}{2} + \Delta\}$, possibly empty. Without loss of generality we assume $\Delta < \frac{1}{4}$ hereafter. In view of the arguments in Robinson (1995), $\widehat{d} \overset{p}{\to} d_0$ if

$$
\sup_{\Theta_1} |T(d)| \overset{p}{\to} 0 \quad \text{and} \quad \Pr\left( \inf_{\Theta_2} S(d) \le 0 \right) \to 0
$$

as $n \to \infty$, where

$$
\begin{aligned}
T(d) &= \log \frac{\widehat{G}(d_0)}{G_0} \\
&\quad - \log \frac{\widehat{G}(d)}{G(d)} - \log \left( \frac{1}{m} \sum_1^m j^{2d-2d_0} \Big/ \frac{m^{2d-2d_0}}{2(d-d_0)+1} \right) \\
&\quad + (2d-2d_0) \left[ \frac{1}{m} \sum_1^m \log j - (\log m - 1) \right].
\end{aligned}
$$

Robinson (1995) shows that the fourth term on the right-hand side is $O(\log m/m)$ uniformly in $d \in \Theta_1$ and

$$
(30) \qquad \sup_{\Theta_1} \left| \frac{2(d-d_0)+1}{m} \sum_1^m \left( \frac{j}{m} \right)^{2d-2d_0} - 1 \right| = O\left( \frac{1}{m^{2\Delta}} \right).
$$

Thus, $\sup_{\Theta_1} |T(d)| \xrightarrow{p} 0$ if

$$(31) \qquad \sup_{\Theta_1} \left| \log \frac{\widehat{G}(d_0)}{G_0} - \log \frac{\widehat{G}(d)}{G(d)} \right| \xrightarrow{p} 0.$$

Let

$$A(d) = \frac{2(d-d_0)+1}{m} \sum_1^m \left(\frac{j}{m}\right)^{2d-2d_0} [\lambda_j^{2d_0} I_x(\lambda_j) - G_0],$$

$$B(d) = \frac{2(d-d_0)+1}{m} G_0 \sum_1^m \left(\frac{j}{m}\right)^{2d-2d_0},$$

from which it follows that

$$\frac{\widehat{G}(d) - G(d)}{G(d)} = \frac{A(d)}{B(d)},$$

$$\log \frac{\widehat{G}(d_0)}{G_0} - \log \frac{\widehat{G}(d)}{G(d)} = \log \left(\frac{B(d)}{B(d_0)}\right) + \log \left(\frac{B(d_0) + A(d_0)}{B(d) + A(d)}\right).$$

By the fact that $\Pr(|\log Y| \geq \varepsilon) \leq 2\Pr(|Y-1| \geq \varepsilon/2)$ for any nonnegative random variable $Y$ and $\varepsilon \leq 1$, (31) holds if

$$(32) \qquad \begin{aligned} &\sup_{\Theta_1} \left| \frac{B(d) - B(d_0)}{B(d_0)} \right| \xrightarrow{p} 0 \quad \text{and} \\ &\sup_{\Theta_1} \left| \frac{B(d_0) - B(d) + A(d_0) - A(d)}{B(d) + A(d)} \right| \xrightarrow{p} 0. \end{aligned}$$

For $d_0 \in (\frac{1}{2}, 1)$, from the arguments in Robinson [(1995), page 1636], $\sup_{\Theta_1} |A(d)|$ is bounded by

$$(33) \qquad \begin{aligned} &\sum_{r=1}^{m-1} \left(\frac{r}{m}\right)^{2\Delta} \frac{1}{r^2} \left| \sum_{j=1}^r [\lambda_j^{2d_0} I_x(\lambda_j) - G_0] \right| \\ &\quad + \frac{1}{m} \left| \sum_{j=1}^m [\lambda_j^{2d_0} I_x(\lambda_j) - G_0] \right|. \end{aligned}$$

Define $D_{nj}(d) = (1 - e^{i\lambda_j})^{-1} \lambda_j^d D_n(e^{i\lambda_j}; \theta)$. Then from Lemmas A.2 and A.3 we have

$$(34) \qquad \begin{aligned} D_{nj}(d) &= e^{(\pi/2)di} + O(\lambda_j) + O(j^{d-2}), \\ |D_{nj}(d)|^2 &= 1 + O(\lambda_j^2) + O(j^{d-2}) \end{aligned}$$

uniformly in $j = 1, \ldots, m$. Hereafter let $I_{xj}$ denote $I_x(\lambda_j)$, let $w_{uj}$ denote $w_u(\lambda_j)$, and similarly for other d.f.t.'s and periodograms. Now

$$\lambda_j^{2d_0} I_{xj} - G_0$$

$$(35) \qquad = \lambda_j^{2d_0} I_{xj} - |D_{nj}(d_0)|^2 I_{uj} + [|D_{nj}(d_0)|^2 - f_u(0)/f_u(\lambda_j)] I_{uj}$$

$$+ [I_{uj} - |C(e^{i\lambda_j})|^2 I_{\varepsilon j}] f_u(0)/f_u(\lambda_j) + f_u(0)(2\pi I_{\varepsilon j} - 1).$$

From Lemma A.1(a), the fact that $||A|^2 - |B|^2| \le |A + B||A - B|$ and the Cauchy–Schwarz inequality we have

$$E|\lambda_j^{2d_0} I_{xj} - |D_{nj}(d_0)|^2 I_{uj}|$$

$$(36) \qquad \le \left( E \left| 2D_{nj}(d_0) w_{uj} - \frac{\lambda_j^{d_0}}{1 - e^{i\lambda_j}} \frac{\widetilde{U}_{\lambda_j n}(\theta_0) + e^{i\lambda_j} X_n}{\sqrt{2\pi n}} \right|^2 \right)^{1/2}$$

$$\times \left( E \left| \frac{\lambda_n^{d_0}}{1 - e^{i\lambda_j}} \frac{\widetilde{U}_{\lambda_j n}(\theta_0) + e^{i\lambda_j} X_n}{\sqrt{2\pi n}} \right|^2 \right)^{1/2},$$

with $\theta_0 = 1 - d_0$. From (34), Lemma A.5(a) and $EI_{uj} = O(1)$ [Robinson (1995), page 1637] the right-hand side is $O(j^{d_0 - 1})$, giving

$$E \left\{ \sum_1^{m-1} \left( \frac{r}{m} \right)^{2\Delta} \frac{1}{r^2} \left| \sum_1^r [\lambda_j^{2d_0} I_{xj} - |D_{nj}(d_0)|^2 I_{uj}] \right| \right\}$$

$$= O(m^{d_0 - 1} + m^{-2\Delta} \log m).$$

For any $\eta > 0$, (34) and Assumption 1 imply that $n$ can be chosen so that

$$||D_{nj}(d_0)|^2 - f_u(0)/f_u(\lambda_j)| \le \eta + O(\lambda_j^2) + O(j^{-1/2}), \qquad j = 1, \ldots, m,$$

and from Robinson [(1995), page 1637], we have

$$E|I_{uj} - |C(e^{i\lambda_j})|^2 I_{\varepsilon j}| = O(j^{-1/2}(\log(j+1))^{1/2}), \qquad j = 1, \ldots, m.$$

It follows that

$$\sum_1^m \left( \frac{r}{m} \right)^{2\Delta} \frac{1}{r^2} \sum_1^r |[|D_{nj}(d_0)|^2 - f_u(0)/f_u(\lambda_j)] I_{uj}$$

$$+ [I_{uj} - |C(e^{i\lambda_j})|^2 I_{\varepsilon j}] f_u(0)/f_u(\lambda_j)|$$

$$= O_p(\eta + m^2 n^{-2} + m^{-2\Delta} \log m).$$

Robinson (1995) shows $\sum_1^m (r/m)^{2\Delta} r^{-2} |\sum_1^r (2\pi I_{\varepsilon j} - 1)| \xrightarrow{p} 0$. Using the same technique, we can show that the second term in (33) is $o_p(1)$, giving $\sup_{\Theta_1} |A(d)| \xrightarrow{p} 0$.

For $d_0 = 1$ first observe that

$$A(d) - \frac{X_n^2}{2\pi n} \frac{2(d - d_0) + 1}{m} \sum_1^m \left(\frac{j}{m}\right)^{2d - 2d_0}$$

$$= \frac{2(d - d_0) + 1}{m} \sum_1^m \left(\frac{j}{m}\right)^{2d - 2d_0} \left[\lambda_j^{2d_0} I_{xj} - G_0 - \frac{X_n^2}{2\pi n}\right].$$

From Lemma A.1(b) we have

(37)
$$\lambda_j^{2d_0} I_{xj} - \frac{X_n^2}{2\pi n}$$

$$= \frac{\lambda_j^2}{|1 - e^{i\lambda_j}|^2} I_{uj} + O(\lambda_j^2) \frac{X_n^2}{2\pi n} - \frac{X_n}{\sqrt{2\pi n}} \frac{\lambda_j^2 2 \operatorname{Re}[e^{i\lambda_j} w_{uj}]}{|1 - e^{i\lambda_j}|^2}.$$

The results in Robinson [(1995), page 1637] imply that

(38)    $$E|w_{uj} - C(e^{i\lambda_j}) w_{\varepsilon j}|^2 = O(j^{-1} \log(j + 1)), \qquad j = 1, \ldots, m.$$

Using a similar decomposition as (35) and the results thereafter, with (37) and (38), we obtain

(39)
$$E\left[\sum_{j=1}^r (\lambda_j^{2d_0} I_{xj} - 2\pi f_u(0) I_{\varepsilon j} - X_n^2/(2\pi n))\right]$$

$$= O(r\eta + r^3 n^{-2} + r^{1/2} \log r),$$

for $1 \le r \le m$. In view of (30) it follows that

$$\sup_{\Theta_1} |A(d) - X_n^2/(2\pi n)| = O_p(\eta + m^2 n^{-2} + m^{-2\Delta} \log m) + o_p(1).$$

Finally, observe that (30) gives $\sup_{\Theta_1} |B(d) - G_0| = O(m^{-2\Delta})$ and (32) follows.

Now we consider $\Theta_2 = \{d : \Delta_1 \le d < d_0 - \frac{1}{2} + \Delta\}$. In a way similar to Robinson [(1995), pages 1638 and 1639] we have

$$\Pr\left(\inf_{\Theta_2} S(d) \le 0\right) \le \Pr\left(\frac{1}{m} \sum_1^m (a_j - 1) \lambda_j^{2d_0} I_{xj} \le 0\right),$$

where $p = \exp(m^{-1} \sum_1^m \log j) \sim m/e$ as $m \to \infty$ and

$$a_j = \begin{cases} (j/p)^{2\Delta - 1}, & \text{for } 1 \le j \le p, \\ (j/p)^{-2d_0 - 1}, & \text{for } p < j \le m, \end{cases}$$

$$\sum_1^m a_j = O(m),$$

$$\sum_1^m a_j{}^2 = O(m^{2-4\Delta}),$$

$$\sum_1^m a_j j^{d_0-1} = O(m^{1-2\Delta}\log m + m^{d_0}),$$

$$\sum_1^m a_j j^{-1/2} = O(m^{1-2\Delta}).$$

Applying (35) and (37) and proceeding as above in conjunction with the fact above and $m^{-1}\sum_1^m(a_j-1)(2\pi I_{\varepsilon j}-1) \overset{p}{\to} 0$ [Robinson (1995), page 1639], we obtain

$$\frac{1}{m}\sum_{j=1}^m(a_j-1)\lambda_j^{2d_0}I_{xj} = \left(G_0 + \frac{X_n^2}{2\pi n}I\{d_0=1\}\right)\frac{1}{m}\sum_1^m(a_j-1) + o_p(1).$$

Choose $\Delta < 1/(2e) < 1/4$ with no loss of generality. Then for sufficiently large $m$ we have $m^{-1}\sum_1^m(a_j-1) > \delta > 0$ and hence

$$\Pr\left(\frac{1}{m}\sum_1^m(a_j-1)\lambda_j^{2d_0}I_{xj} \le 0\right) \to 0$$

as $n \to \infty$. Therefore, $\hat{d} \overset{p}{\to} d_0$, giving the stated result.

For the limit of $\widehat{G}(d)$, recall $\widehat{G}(d) = G(d) + A(d)G(d)/B(d)$, $\hat{d} \overset{p}{\to} d_0$, $G(\hat{d}) \overset{p}{\to} G_0$ and $B(\hat{d}) \overset{p}{\to} G_0$. The required result follows because

$$\sup_{\Theta_1}\left|A(d) - \frac{X_n^2}{2\pi n}I\{d_0=1\}\right| \overset{p}{\to} 0,$$

and $X_n^2/(2\pi n) = G_0(X_n^\varepsilon)^2 + o_p(1) \overset{d}{\to} G_0\chi_1^2$ from a standard martingale CLT. $\square$

PROOF OF THEOREM 3.2. Define $G(d) = G_0 m^{-1}\sum_1^m \lambda_j^{2d-2}$ and $S(d) = R(d) - R(1)$. For $0 < \Delta < \frac{1}{4}$ define $\Theta_1 = \{d : \frac{1}{2} + \Delta \le d \le \Delta_2\}$ and $\Theta_2 = \{d : \Delta_1 \le d < \frac{1}{2} + \Delta\}$, possibly empty. Then, by the same line of arguments as above, $\hat{d} \overset{p}{\to} 1$ if

$$\sup_{\Theta_1}|T(d)| \overset{p}{\to} 0 \quad \text{and} \quad \Pr\left(\inf_{\Theta_2} S(d) \le 0\right) \to 0$$

as $n \to \infty$, where

$$T(d) = \log\frac{\widehat{G}(1)}{G_0} - \log\frac{\widehat{G}(d)}{G(d)} - \log\left(\frac{1}{m}\sum_1^m j^{2d-2} \Big/ \frac{m^{2d-2}}{2(d-1)+1}\right)$$

$$+ (2d-2)\left[\frac{1}{m}\sum_1^m \log j - (\log m - 1)\right];$$

$\sup_{\Theta_1} |T(d)| \overset{p}{\to} 0$ if

(40) $$\sup_{\Theta_1} \left| \log \frac{\widehat{G}(1)}{G_0} - \log \frac{\widehat{G}(d)}{G(d)} \right| \overset{p}{\to} 0.$$

Let

$$A(d) = \frac{2d-1}{m} \sum_1^m \left(\frac{j}{m}\right)^{2d-2} j^{2-2d_0} \lambda_j^{2d_0} I_{xj},$$

$$B(d) = \frac{2d-1}{m} G_0 \sum_1^m \left(\frac{j}{m}\right)^{2d-2}.$$

Then a little algebra shows $[\widehat{G}(d)/G(d)] = (2\pi/n)^{2-2d_0}[A(d)/B(d)]$, giving

$$\log \frac{\widehat{G}(1)}{G_0} - \log \frac{\widehat{G}(d)}{G(d)} = \log\left(\frac{B(d)}{G_0}\right) - \log\left(\frac{A(d)}{A(1)}\right).$$

Therefore, (40) holds if

(41) $$\sup_{\Theta_1} \left| \frac{A(d) - A(1)}{A(1)} \right| \overset{p}{\to} 0 \quad \text{and} \quad \sup_{\Theta_1} \left| \frac{B(d) - G_0}{G_0} \right| \overset{p}{\to} 0.$$

We proceed to approximate $A(d)$ by

$$A_1(d) = \frac{2d-1}{m} \sum_1^m \left(\frac{j}{m}\right)^{2d-2} j^{2-2d_0} \lambda_j^{2d_0} |1 - e^{i\lambda_j}|^{-2} \frac{X_n^2}{2\pi n}$$

$$= (2\pi)^{2d_0-3} C(1)^2 n^{1-2d_0} (X_n^\varepsilon)^2 + o_p(1) \qquad \text{uniformly in } d \in \Theta_1,$$

where the second equality follows from

$$j^{2-2d_0} \lambda_j^{2d_0} |1 - e^{i\lambda_j}|^{-2} = (2\pi)^{2d_0-2} n^{2-2d_0} (1 + O(\lambda_j^2)),$$

(30) and Lemma A.5(a2). For $d_0 \in (1, \frac{3}{2}]$, from Lemmas A.1(a) and A.5, similarly as in (36) we obtain [$D_{nj}(d)$ is defined in the proof of Theorem 3.1]

$$E\left| j^{2-2d_0} \lambda_j^{2d_0} I_{xj} - j^{2-2d_0} \lambda_j^{2d_0} |1 - e^{i\lambda_j}|^{-2} \frac{X_n^2}{2\pi n} \right|$$

$$\leq \left( E\left| j^{1-d_0} D_{nj}(d_0) w_{uj} - \frac{j^{1-d_0} \lambda_j^{d_0}}{1 - e^{i\lambda_j}} \frac{\widetilde{U}_{\lambda_j n}(\theta_0) + 2e^{i\lambda_j} X_n}{\sqrt{2\pi n}} \right|^2 \right)^{1/2}$$

$$\times \left( E\left| j^{1-d_0} D_{nj}(d_0) w_{uj} - \frac{j^{1-d_0} \lambda_j^{d_0}}{1 - e^{i\lambda_j}} \frac{\widetilde{U}_{\lambda_j n}(\theta_0)}{\sqrt{2\pi n}} \right|^2 \right)^{1/2}$$

$$= O(j^{1-d_0}).$$

It follows that $\sup_{\Theta_1} |A(d) - A_1(d)| = O_p(m^{1-d_0} + m^{-2\Delta})$, and uniformly in $\Theta_1$ we have

$$\frac{A(d) - A(1)}{A(1)} = \frac{o_p(1)}{(2\pi)^{2d_0-3}C(1)^2 n^{1-2d_0}(X_n^\varepsilon)^2 + o_p(1)}$$

$$= \frac{o_p(1)}{(2\pi)^{2d_0-3}C(1)^2(\sum_1^n y_t)^2 + o_p(1)},$$

where $y_t = n^{1/2-d_0}(d_0)_{n-t}\varepsilon_t/(n-t)!$. Assumption 3 implies

$$\sum_1^n E(y_t^2|F_{t-1}) \to \Phi_1 = \Gamma(d_0)^{-2}(2d_0-1)^{-1},$$

$$\sum_1^n E(y_t^2 I\{|y_t| > \delta\}) \to 0 \qquad \text{for all } \delta > 0.$$

Therefore, from a standard martingale CLT we have $n^{1/2-d_0}X_n^\varepsilon \xrightarrow{d} N(0, \Phi_1)$. Thus,

$$\sup_{\Theta_1} \left| \frac{A(d) - A(1)}{A(1)} \right| \xrightarrow{p} 0,$$

and $\sup_{\Theta_1} |[B(d) - G_0]/G_0| \to 0$ as before, thereby establishing (41).

Next, consider $\Theta_2 = \{d : \Delta_1 \le d < \frac{1}{2} + \Delta\}$. Let $p = \exp(m^{-1}\sum_1^m \log j)$. Then $S(d) = \log\{\widehat{D}(d)/\widehat{D}(1)\}$, where $\widehat{\widehat{D}}(d) = m^{-1}\sum_1^m (j/p)^{2d-2}j^2 I_{xj}$. It follows that

$$\inf_{\Theta_2} \widehat{D}(d) \ge \frac{1}{m}\sum_1^m a_j j^2 I_{xj},$$

where

$$a_j = \begin{cases} (j/p)^{2\Delta-1}, & \text{for } 1 \le j \le p, \\ (j/p)^{-3}, & \text{for } p < j \le m. \end{cases}$$

Then

$$(42) \qquad \Pr\left(\inf_{\Theta_2} S(d) \le 0\right) \le \Pr\left(\frac{1}{m}\sum_1^m (a_j-1)j^{2-2d_0}\lambda_j^{2d_0} I_{xj} \le 0\right).$$

In view of the fact that $\sum_1^m a_j = O(m)$, $\sum_1^m a_j j^{1-d_0} = O(m^{1-2\Delta}\log m + m^{-d_0})$, $\sum_1^m a_j j^{-1/2} = O(m^{1-2\Delta})$, we obtain similarly to before

$$\frac{1}{m}\sum_1^m (a_j-1)j^{2-2d_0}\lambda_j^{2d_0} I_{xj}$$

$$= \frac{1}{m}\sum_1^m (a_j-1)j^{2-2d_0}\lambda_j^{2d_0}|1 - e^{i\lambda_j}|^{-2}\frac{X_n^2}{2\pi n} + o_p(1)$$

$$= (2\pi)^{2d_0-3} n^{1-2d_0} X_n^2 \frac{1}{m} \sum_1^m (a_j - 1) + o_p(1).$$

Since $m^{-1} \sum_1^m (a_j - 1) > \delta > 0$ for sufficiently large $m$ by choosing $\Delta < 1/(2e)$, we obtain $\Pr(m^{-1} \sum_1^m (a_j - 1) j^{2-2d_0} \lambda_j^{2d_0} I_{xj} \le 0) \to 0$ as $n \to \infty$ and hence $\hat{d} \xrightarrow{p} 1$, giving the stated result. For $d_0 \in (\frac{3}{2}, \frac{5}{2}]$, from Lemma A.1(b) we have

$$(1 - e^{i\lambda_j}) w_{xj} = w_{\Delta x_j} - n^{d-1} n^{1/2-d} X_n e^{i\lambda_j}/\sqrt{2\pi}.$$

Because $E|w_{\Delta x_j}|^2 = O(n^{2d-2} s^{-1})$ from $\Delta X_t \sim I(d-1)$ and $n^{1/2-d} X_n$ converges to a Gaussian random variable, the stochastic behavior of $w_{xj}$ is dominated by $X_n$. Hence, the required result follows from the same line of argument as above, and the results for larger $d_0$ are derived similarly. $\square$

PROOF OF THEOREM 4.1. We follow the same line of approach as the proof of Theorem 2 of Robinson (1995). Theorem 3.1 holds under the current conditions and implies that with probability approaching 1, as $n \to \infty$, $\hat{d}$ satisfies

$$(43) \qquad 0 = R'(\hat{d}) = R'(d_0) + R''(d^*)(\hat{d} - d_0),$$

where $|d^* - d_0| \le |\hat{d} - d_0|$. Now

$$R''(d) = \frac{4[\hat{F}_2(d)\hat{F}_0(d) - \hat{F}_1^2(d)]}{\hat{F}_0^2(d)} = \frac{4[\hat{E}_2(d)\hat{E}_0(d) - \hat{E}_1^2(d)]}{\hat{E}_0^2(d)},$$

$$\hat{F}_k(d) = \frac{1}{m} \sum_1^m (\log j)^k \lambda_j^{2d} I_{xj}, \qquad \hat{E}_k(d) = \frac{1}{m} \sum_1^m (\log j)^k j^{2d} I_{xj}.$$

As pointed out by Andrews and Sun [(2001), page 21], since $\hat{F}_k(d_0) = O_p((\log m)^k)$ as shown below, we need to show $\hat{E}_k(d^*) - \hat{E}_k(d_0) = o_p(n^{2d_0}(\log m)^{-k})$ rather than $o_p(n^{2d_0})$ as in (4.4) in Robinson. Fix $\varepsilon > 0$ and choose $n$ so that $2\varepsilon < (\log m)^2$. Let $M = \{d : (\log m)^6 |d - d_0| \le \varepsilon\}$. As in Robinson [(1995), page 1642] we have

$$(44) \qquad \Pr(|\hat{E}_k(d^*) - \hat{E}_k(d_0)| > (2\pi/n)^{-2d_0}(\log m)^{-k})$$
$$\le \Pr(\hat{G}(d_0) > (\log m)^{5-2k}/(2e\varepsilon)) + \Pr((\log m)^6 |d^* - d_0| > \varepsilon).$$

The first probability tends to 0 because $\hat{G}(d_0) \xrightarrow{p} G_0$. The second probability tends to 0 if

$$(45) \qquad \sup_{\Theta_1} \left| \frac{B(d) - B(d_0)}{B(d_0)} \right| + \sup_{\Theta_1} \left| \frac{B(d_0) - B(d) + A(d_0) - A(d)}{B(d) + A(d)} \right|$$
$$= o_p((\log m)^{-12}).$$

Using (35) and Assumption 1', we obtain, for $1 \le r \le m$,

$$E\left\{\sum_1^r (\lambda_j^{2d_0} I_{xj} - f_u(0) 2\pi I_{\varepsilon j})\right\} = O(r^{d_0} + r^{\beta+1} n^{-\beta}),$$

and Robinson [(1995), (4.9)] shows that $\sum_1^r (2\pi I_{\varepsilon j} - 1) = O_p(r^{1/2})$. In conjunction with $\sup_{\Theta_1} |B(d) - G_0| = O(m^{-2\Delta})$ they give (45). It follows that

$$R''(d^*) = 4[\widehat{F}_2(d_0)\widehat{F}_0(d_0) - \widehat{F}_1^2(d_0)][\widehat{F}_0^2(d_0)]^{-1} + o_p(1) = 4 + o_p(1),$$

where the second equality follows from $\widehat{F}_k(d_0) = G_0 m^{-1}\sum_1^m (\log j)^k + o_p((\log m)^{-3})$, obtained similarly as $A(d_0) \xrightarrow{p} 0$. Next we consider the first term on the right-hand side of (43). Now

$$m^{1/2} R'(d_0) = 2m^{-1/2}\sum_1^m \nu_j [\lambda_j^{2d_0} I_{xj} - G_0][G_0 + o_p(1)]^{-1},$$

where $\nu_j = \log j - m^{-1}\sum_1^m \log j$ and $\sum_1^m \nu_j = 0$. From Lemmas A.1 and A.5, (34) and (38), we have

$$\sum_1^m \nu_j \lambda_j^{2d_0} I_{xj} = \sum_1^m \nu_j I_{uj} + (2\pi n)^{-1} X_n^2 \sum_1^m \nu_j \lambda_j^{2d_0} |1 - e^{i\lambda_j}|^{-2} + 2\operatorname{Re}[T_n] + R_n,$$

where

$$T_n = \sum_1^m \nu_j D_{nj}^*(d_0) C^*(e^{i\lambda_j}) w_{\varepsilon j}^* \lambda_j^{d_0} (1 - e^{i\lambda_j})^{-1} C(1)\widetilde{\varepsilon}_{\lambda_j n}(\theta_0)(2\pi n)^{-1/2}$$

and

$$R_n = O_p((\log m)(m^{d_0-1/2}\log m + m^{1/2}(\log n)^{-2} + (\log m)^2 + m^3 n^{-2}))$$

$$= o_p(m^{1/2}).$$

Rewrite $T_n/C(1)$ as

$$\sum_1^m \nu_j D_{nj}^*(d_0) C^*(e^{i\lambda_j}) \lambda_j^{d_0}(1 - e^{i\lambda_j})^{-1}(2\pi n)^{-1}$$
$$\times \sum_{p=0}^{n-1} \widetilde{\theta}_{\lambda_j p} e^{-ip\lambda_j}\varepsilon_{n-p}\sum_{q=0}^{n-1} e^{iq\lambda_j}\varepsilon_{n-q}.$$
(46)

Since $\varepsilon_t$ is a martingale difference sequence, $E|T_n|^2$ is bounded by

(47) $$\frac{1}{n^2}\sum_{j=1}^m \sum_{k=1}^m |\nu_j||\nu_k|\lambda_j^{d_0-1}\lambda_k^{d_0-1}\sum_{p=0}^{n-1} |\widetilde{\theta}_{\lambda_j p}||\widetilde{\theta}_{-\lambda_k p}|$$

$$(48) \qquad + \frac{1}{n^2} \sum_{j=1}^{m} |\nu_j| \lambda_j^{d_0-1} \sum_{p=0}^{n-1} |\widetilde{\theta}_{\lambda_j p}| \sum_{k=1}^{m} |\nu_k| \lambda_k^{d_0-1} \sum_{q=0}^{n-1} |\widetilde{\theta}_{-\lambda_k q}|$$

$$+ \frac{1}{n^2} \sum_{j=1}^{m} \sum_{k=1}^{m} |\nu_j||\nu_k| \lambda_j^{d_0-1} \lambda_k^{d_0-1}$$

$$(49) \qquad \times \sum_{p=0, p \neq q}^{n-1} |\widetilde{\theta}_{\lambda_j p}||\widetilde{\theta}_{-\lambda_k p}| \left| \sum_{q=0}^{n-1} e^{iq(\lambda_j - \lambda_k)} \right|,$$

and (47) and (48) are bounded by, respectively,

$$(\log m)^2 \sum_{1}^{m} \sum_{1}^{m} j^{d_0-1} k^{d_0-1} n^{-2d_0} \sum_{0}^{n-1} |p|_+^{2d_0-2} = O(n^{-1} m^{2d_0} (\log m)^2),$$

$$\left[ \log m \sum_{1}^{m} j^{d_0-1} n^{-d_0} \left( \sum_{0}^{n/j} |p|_+^{d_0-1} + \sum_{n/j}^{n-1} p^{d_0-2} n j^{-1} \right) \right]^2 = O((\log m)^4).$$

In view of the fact that $\sum_{0}^{n-1} e^{iq(\lambda_j - \lambda_k)} = nI\{j = k\}$, (49) is bounded by

$$(\log m)^2 n^{-1} \sum_{1}^{m} \lambda_j^{2d_0-2} \sum_{0}^{n-1} |p|_+^{2d_0-2} = O(m^{2d_0-1} (\log m)^2),$$

giving $T_n = O_p(n^{-1/2} m^{d_0} \log m + (\log m)^2 + m^{d_0-1/2} \log m) = o_p(m^{1/2})$.

From Lemma A.5 and the fact that $\sum_{1}^{m} \nu_j j^{2d_0-2} = (2d_0 - 1)^{-2}(2d_0 - 2)m^{2d_0-1} + O(\log m)$, we obtain

$$(2\pi n)^{-1} X_n^2 \sum_{1}^{m} \nu_j \lambda_j^{2d_0} |1 - e^{i\lambda_j}|^{-2} = \Xi G_0 m^{2d_0-1} [n^{1-2d_0} (X_n^\varepsilon)^2 + o_p(1)],$$

where $\Xi = (2\pi)^{2d_0-2}(2d_0-1)^{-2}(2d_0-2)$. Robinson [(1995), page 1644] shows that $\sum_{1}^{m} \nu_j I_{uj} = G_0 \sum_{1}^{m} \nu_j I_{\varepsilon j} + o_p(m^{1/2})$. Therefore,

$$m^{1/2} R'(d_0) = 2m^{-1/2} \sum_{1}^{m} \nu_j [2\pi I_{\varepsilon j} - 1]$$

$$- 2\Xi m^{2d_0-3/2} [n^{1-2d_0} (X_n^\varepsilon)^2 + o_p(1)] + o_p(1).$$

The first term on the right-hand side converges to a $N(0,4)$ random variable by Robinson (1995). For $d_0 \in (\frac{1}{2}, \frac{3}{4})$, the second term on the right-hand side is $o_p(1)$, and the required result follows. For $d_0 \in (\frac{3}{4}, 1)$, we have

$$m^{2-2d_0} R'(d_0) = 2\Xi n^{1-2d_0} (X_n^\varepsilon)^2 + o_p(1) = 2\Xi \left( \sum_{1}^{n} y_t \right)^2 + o_p(1),$$

where $y_t = n^{1/2-d_0}(d_0)_{n-t} \varepsilon_t / (n-t)!$, suppressing reference to $n$ in $y_t$. Since $\sum_{1}^{n} E(y_t^2 | F_{t-1}) \to \Phi_1 = \Gamma(d_0)^{-2}(2d_0 - 1)^{-1}$ and $\sum_{1}^{n} E y_t^4 = O(n^{-1})$, from a standard martingale CLT we obtain $\sum_1^n y_t \xrightarrow{d} N(0, \Phi_1)$, giving $m^{2-2d_0} R'(d_0) \xrightarrow{d} 2\Xi\Phi_1\chi_1^2$ and the required result. When $d_0 = \frac{3}{4}$, $m^{1/2} R'(d_0) = 2\sum_1^n z_t + 2\Xi(\sum_1^n y_t)^2 + o_p(1)$, where $y_t$ is defined above, $z_1 = 0$ and, for $t \geq 2$,

$$z_t = \varepsilon_t \sum_{s=1}^{t-1} \varepsilon_s c_{t-s}, \qquad c_s = 2n^{-1}m^{-1/2} \sum_1^m \nu_j \cos(s\lambda_j),$$

and $\xi_t = (z_t, y_t)'$ form a zero-mean martingale difference array; hence $\sum_1^n \xi_t \xrightarrow{d} N(0, \text{diag}(1, \Phi_1))$ if, for any nonrandom $(2 \times 1)$ vector $\alpha$,

$$(50) \qquad \sum_1^n E[(\alpha'\xi_t)^2 | F_{t-1}] - \alpha' \text{diag}(1, \Phi_1)\alpha \xrightarrow{p} 0,$$

$$(51) \qquad \sum_1^n E[(\alpha'\xi_t)^2 I(|\alpha'\xi_t| > \delta)] \to 0 \qquad \text{for all } \delta > 0.$$

Robinson shows $\sum_1^n E(z_t^2 | F_{t-1}) - 1 \xrightarrow{p} 0$ and $\sum_1^n E z_t^4 \to 0$. In conjunction with $\sum_1^n E y_t^4 \to 0$, (51) is satisfied. Since $\sum_1^n E(y_t^2 | F_{t-1}) \to \Phi_1$, (50) holds if

$$\sum_1^n E[y_t z_t | F_{t-1}] = n^{1/2-d_0} \sum_{t=2}^n \frac{(d_0)_{n-t}}{(n-t)!} \sum_{s=1}^{t-1} \varepsilon_s c_{t-s} \xrightarrow{p} 0.$$

The term in the middle has mean zero and variance bounded by

$$n^{1-2d_0} \sum_{t=2}^n |n-t|_+^{d_0-1} \sum_{u=2}^n |n-u|_+^{d_0-1} \sum_{s=1}^{\min\{t-1,u-1\}} c_{t-s} c_{u-s}$$

$$= n^{1-2d_0} \sum_2^n |n-t|_+^{2d_0-2} \sum_1^{t-1} c_{t-s}^2$$

$$+ 2n^{1-2d_0} \sum_3^n |n-t|_+^{d_0-1} \sum_2^{t-1} (n-u)^{d_0-1} \sum_1^{u-1} c_{t-r} c_{u-r}.$$

Robinson [1995, page 1646] shows $c_s = c_{n-s}$, $|c_s| = O(m^{-1/2}s^{-1} \log m)$ for $1 \leq s \leq n/2$, $|c_s| = O(m^{1/2}n^{-1} \log m)$, and $\sum_1^n c_s^2 = O(n^{-1}(\log m)^2)$. Therefore, the first term on the right-hand side is bounded by $n^{-1}(\log m)^2$, and the second term on the right-hand side is bounded by

$$\left(\sum_1^n c_t^2\right)^{1/2} n^{1-2d_0} \sum_3^n |n-t|_+^{d_0-1} \sum_2^{t-1} (n-u)^{d_0-1} \left(\sum_{t-u+1}^{t-1} c_r^2\right)^{1/2}$$

$$= O\left((\log m)^2 m^{-1/2} n^{1/2-2d_0}\right)$$

$$\times \sum_3^n |n-t|_+^{d_0-1} \sum_2^{t-1} (n-u)^{d_0-1}((t-u)^{-1/2} + |n-t|_+^{-1/2})\Big)$$

$$= O(m^{-1/2}(\log m)^2),$$

giving (50), thereby completing the proof of the theorem. $\square$

PROOF OF THEOREM 4.2. We follow the approach and notation of the proof of Theorem 4.1. First, from Assumption 1′ (39) is strengthened to

$$(52) \qquad E\left[\sum_1^r (\lambda_j^{2d_0} I_{xj} - 2\pi f_u(0) I_{\varepsilon j} - X_n^2/2\pi n)\right]$$
$$= O(r^{\beta+1} n^{-\beta} + r^{1/2} \log r),$$

for $1 \le r \le m$. It follows that

$$\sup_{\Theta_1} |A(d) - X_n^2/2\pi n| = O_p(m^\beta n^{-\beta} + m^{-2\Delta} \log m),$$

and thus (45) holds. Since $\widehat{G}(d_0) = G_0 + X_n^2/(2\pi n) + o_p(1)$, the first probability in (44) tends to 0 and $R''(d^*) = 4[\widehat{F}_2(d_0)\widehat{F}_0(d_0) - \widehat{F}_1^2(d_0)][\widehat{F}_0(d_0)]^{-2} + o_p(1)$. Using (52) and the fact that $\sum_1^r (2\pi I_{\varepsilon j} - 1) = O_p(r^{1/2})$, we obtain

$$\widehat{F}_k(d_0) - \left(G_0 + \frac{X_n^2}{2\pi n}\right)\left[\frac{1}{m}\sum_1^m (\log j)^k\right]$$
$$= O_p(m^\beta n^{-\beta}(\log m)^2 + m^{-1/2}(\log m)^3),$$

giving $R''(d^*) \xrightarrow{p} 4$. Now

$$m^{1/2} R'(d_0) = \frac{2m^{-1/2}\sum_1^m \nu_j \lambda_j^{2d_0} I_{xj}}{\widehat{G}(d_0)}$$
$$= \frac{2m^{-1/2}\sum_1^m \nu_j \lambda_j^{2d_0} I_{xj}}{G_0(1 + n^{-1}(X_n^\varepsilon)^2) + o_p(1)}.$$

The numerator is equal to

$$2m^{-1/2}\sum_1^m \nu_j(1 + O(\lambda_j^2))|w_{uj} - e^{i\lambda_j} X_n/\sqrt{2\pi n}|^2$$

$$= 2m^{-1/2}\sum_1^m \nu_j I_{uj} - 2m^{-1/2}(X_n/\sqrt{2\pi n})2\left[\sum_1^m \nu_j e^{-i\lambda_j} w_{uj}\right]$$

$$+ O_p(m^{5/2} n^{-2} \log m)$$

$$= 2m^{-1/2} G_0 \sum_1^m \nu_j [2\pi I_{\varepsilon j} - 1]$$

$$- 2m^{-1/2}G_0 n^{-1/2} X_n^{\varepsilon} \sum_1^m \nu_j 2\operatorname{Re}[\sqrt{2\pi}\,w_{\varepsilon j}] + o_p(1),$$

where the third line follows from Robinson [(1995), (4.8)], Lemma A.6 and Assumptions 1′ and 6′. It follows that

$$m^{1/2}R'(d_0) = \frac{2\sum_1^n z_t - 2\sum_1^n y_t \sum_1^n x_t + o_p(1)}{1 + (\sum_1^n y_t)^2 + o_p(1)},$$

where $y_t = n^{-1/2}\varepsilon_t$, $x_t = n^{1/2}\varepsilon_t c_t$ and $z_t$ and $c_t$ are defined in the proof of Theorem 4.1. Therefore, $W_n = \sum_1^n (z_t, y_t, x_t)' \xrightarrow{d} W \sim N(0, \operatorname{diag}(1,1,2))$ if

$$(53) \qquad \sum_1^n E[(z_t, y_t, x_t)'(z_t, y_t, x_t)|F_{t-1}] \xrightarrow{p} \operatorname{diag}(1,1,2),$$

$$(54) \qquad \sum_1^n Ez_t^4 + \sum_1^n Ey_t^4 + \sum_1^n Ex_t^4 \to 0.$$

We have already shown $\sum_1^n Ez_t^4 + \sum_1^n Ey_t^4 \to 0$ in the proof of Theorem 4.1. Since

$$\sum_1^n Ex_t^4 = n^2 \sum_1^n c_t^4$$

$$= O\left(n^{-2}\sum_1^{n/m} m^2(\log m)^4 + n^2 \sum_{n/m}^{\infty} m^{-2}s^{-4}(\log m)^4\right)$$

$$= O(n^{-1}m(\log m)^4),$$

(54) holds. To show (53), $\sum_1^n E[(z_t, y_t)'(z_t, y_t)|F_{t-1}] \xrightarrow{p} \operatorname{diag}(1,1)$ has already been shown above, and

$$\sum_1^n E(x_t^2|F_{t-1})$$

$$= 4n^{-1}m^{-1}\sum_1^n \left(\sum_1^m \nu_j \cos(s\lambda_j)\right)^2$$

$$= 4n^{-1}m^{-1}\sum_1^m \nu_j^2 \sum_1^n \cos^2(s\lambda_j)$$

$$\qquad + 2n^{-1}m^{-1}\sum_{j \neq k}\sum \nu_j \nu_k \sum_1^n [\cos\{s(\lambda_j + \lambda_k)\} + \cos\{s(\lambda_j - \lambda_k)\}]$$

$$= 2n^{-1}m^{-1}\sum_1^m \nu_j^2 \sum_1^n [1 + \cos(2s\lambda_j)]$$

$\to 2,$

since $\sum_1^m \nu_j^2 \sim m$. Furthermore, $\sum_1^n E(x_t y_t | F_{t-1}) = \sum_1^n c_t = 0$ and

$$\sum_1^n E(x_t z_t | F_{t-1}) = \sum_{t=2}^n n^{1/2} c_t \sum_{s=1}^{t-1} \varepsilon_s c_{t-s} \xrightarrow{p} 0,$$

because the right-hand side has mean zero and variance

$$n \sum_{t=2}^n c_t \sum_{u=2}^n c_u \sum_{s=1}^{\min\{t-1,u-1\}} c_{t-s} c_{u-s}$$

$$= n \sum_2^n c_t^2 \sum_1^{t-1} c_{t-s}^2 + 2n \sum_3^n c_t \sum_2^{t-1} c_u \sum_1^{u-1} c_{t-s} c_{u-s}.$$

The first term is $O(n^{-1}(\log m)^4)$, and the second term is bounded by

$$n \left( \sum_1^n |c_s| \right)^2 \left( \sum_1^n c_t^2 \right)$$

$$= O\left( (\log m)^2 \left( m^{-1/2} \log m + m^{-1/2} \log m \sum_{n/m}^{n/2} s^{-1} \right)^2 \right)$$

$$= O(m^{-1}(\log m)^6).$$

Thus (53) holds and $W_n \xrightarrow{d} W$. Therefore, from the continuous mapping theorem

$$m^{1/2} R'(d_0) \xrightarrow{d} \frac{2W_1 - 2W_2 W_3}{1 + (W_2)^2} \sim N\left( 0, \frac{4[1 + 2(W_2)^2]}{[1 + (W_2)^2]^2} \right)$$

conditional on $W_2$, and unconditionally

$$m^{1/2}(\hat{d} - d_0) \xrightarrow{d} \int_{-\infty}^{\infty} N(0, \tfrac{1}{4}(1 + 2h^2)(1 + h^2)^{-2}) \phi(h) \, dh,$$

where $\phi(\cdot)$ is the standard normal p.d.f., giving the stated result. $\square$

PROOF OF THEOREM 5.1. The argument follows the approach of the proof of Theorem 3.2. First we consider the case $\alpha \le d_0 - \frac{1}{2}$. Since $\alpha > \frac{1}{2}$, $d_0 > 1$ must hold. Let

$$A_1(d) = \frac{2d - 1}{m} \sum_1^m \left( \frac{j}{m} \right)^{2d-2} \frac{j^{2-2d_0} \lambda_j^{2d_0}}{|1 - e^{i\lambda_j}|^2} \frac{1}{2\pi n} |\mu n^\alpha + e^{i\lambda_j} X_n^0|^2$$

$$= (2\pi)^{2d_0 - 3} [\mu n^{\alpha - d_0 + 1/2} + n^{1/2 - d_0} X_n^0]^2$$

$$\quad + o_p(1) \qquad \text{(uniformly for } \Theta_1)$$

$$\xrightarrow{d} [C_1 n^{\alpha - d_0 + 1/2} + C_2 N(0,1)]^2,$$

for generic nonzero constants $C_1$ and $C_2$, where the second equality and convergence in distribution follow by the same argument as before. Define the other quantities as in the proof of Theorem 3.2. Because

$$E|j^{2-2d_0}\lambda_j^{2d_0}I_{xj} - j^{2-2d_0}\lambda_j^{2d_0}|1 - e^{i\lambda_j}|^{-2}(2\pi n)^{-1}|\mu n^\alpha + e^{i\lambda_j}X_n^0|^2|$$
$$= O(j^{1-d_0} + j^{-1/2} + n^{\alpha - d_0 - 1/2}j),$$

$\sup_{\Theta_1}|A(d) - A_1(d)| \overset{p}{\to} 0$ follows, giving (41); $\Pr(\inf_{\Theta_2} S(d) \leq 0) \to 0$ is obtained similarly and we establish $\hat{d} \overset{p}{\to} 1$.

Next consider the case $\alpha > d_0 - \frac{1}{2}$. Define

$$A(d) = \frac{2d-1}{m}\sum_1^m \left(\frac{j}{m}\right)^{2d-2}n^{1-2\alpha}\lambda_j^2 I_{xj},$$

$$B(d) = \frac{2d-1}{m}G_0\sum_1^m \left(\frac{j}{m}\right)^{2d-2},$$

$$A_1(d) = \frac{2d-1}{m}\sum_1^m \left(\frac{j}{m}\right)^{2d-2}\mu^2\frac{\lambda_j^2}{|1 - e^{i\lambda_j}|^2}\frac{1}{2\pi}$$

$$\to (2\pi)^{-1}\mu^2 \qquad \text{uniformly for } \Theta_1,$$

and define the other quantities as in the proof of Theorem 3.2. Then it follows that

$$\widehat{G}(d)/G(d) = n^{2\alpha-1}A(d)/B(d),$$

and $\sup_{\Theta_1}|T(d)| \overset{p}{\to} 0$ follows if $\sup_{\Theta_1}|[A(d) - A(1)]/A(1)| \overset{p}{\to} 0$. Since

$$E|n^{1-2\alpha}\lambda_j^2 I_{xj} - \lambda_j^2|1 - e^{i\lambda_j}|^{-2}(2\pi)^{-1}\mu^2| = O(j^{-1/2} + jn^{-1} + n^{d_0-1/2-\alpha}),$$

$\sup_{\Theta_1}|A(d) - A_1(d)| \overset{p}{\to} 0$ follows, giving (41); $\Pr(\inf_{\Theta_2} S(d) \leq 0) \to 0$ is obtained similarly, and we establish $\hat{d} \overset{p}{\to} 1$.   □

**Acknowledgments.** The authors thank the Editor, an Associate Editor, two referees, Donald Andrews and Marcus Chambers for helpful comments on earlier versions. In particular, comments by a referee led to a substantial improvement of the results in the paper. Simulations and empirical computations were performed in MATLAB.

COWLES FOUNDATION
FOR RESEARCH IN ECONOMICS
YALE UNIVERSITY
P. O. BOX 208281
NEW HAVEN, CONNECTICUT 06520-8281
USA
E-MAIL: peter.phillips@yale.edu

DEPARTMENT OF ECONOMICS
QUEEN'S UNIVERSITY
KINGSTON, ONTARIO
CANADA K7L 3N6
E-MAIL: shimotsu@qed.econ.queensu.ca